%% file: Triple_dqds1.tex
\documentclass[final]{siamltex}

\usepackage{amsfonts}
\usepackage{amscd}
\usepackage{amsmath}
\usepackage{amssymb}
\usepackage{epic}
\usepackage{multirow}
\usepackage{makeidx}
\usepackage{mathscinet}
\usepackage{hyperref}
\usepackage{graphicx}
\usepackage{subfigure}
\usepackage{fancyhdr}
\usepackage[english]{babel}
\usepackage[latin1]{inputenc}
\usepackage{srcltx}
\usepackage{pmat}

\newcommand{\bs}{\boldsymbol}

\newcommand{\matR}{\mathbb{R}}

\DeclareMathOperator{\cond}{cond}
\DeclareMathOperator{\tridiag}{tridiag}

\DeclareMathOperator{\trace}{trace}

\DeclareMathOperator{\offdiag}{offdiag}
\DeclareMathOperator{\spectrum}{spectrum}

\def\ds{\displaystyle}

\def\norm#1{\|#1\|}

\title{Real \textsl{dqds} for the nonsymmetric tridiagonal eigenvalue problem
\thanks{
The first author is supported by
FEDER Funds through ``Programa Operacional Factores de Competitividade - COMPETE" and by Portuguese Funds through FCT - "Funda\c{c}\~{a}o para a Ci\^{e}ncia e a Tecnologia", within the Project PEst-CMAT/UI0013/2011.
}}

\author{Carla Ferreira\thanks{Mathematics and Applications Department, University of Minho, 4710-057 Braga
        ({\tt caferrei@math.uminho.pt})}
        \and
        Beresford Parlett\thanks{Department of Mathematics and the Computer Science Division of the Electrical Engineering and Computer Science Department, University of California, Berkeley, California 94720 ({\tt parlett@math.berkeley.edu})}}

\begin{document}

\maketitle

\begin{abstract}
We present a new transform \textsl{triple dqds} to help to compute the eigenvalues of a real tridiagonal
matrix $C$ using real arithmetic. The algorithm uses the real \textsl{dqds} transform to shift by a real
number and \textsl{tridqds} to shift by a complex conjugate pair. We present what seems to be a new criteria for
splitting the current pair $L,U$. The algorithm rejects any transform which suffers from excessive element growth
and then tries a new transform. Our numerical tests show that the algorithm is about 100 times faster than the Ehrlich-Aberth
method of D. A. Bini, L. Gemignani and F. Tisseur. Our code is comparable in performance to a complex \textsl{dqds} code and is
sometimes 3 times faster.
\end{abstract}



\pagestyle{myheadings}
\thispagestyle{plain}
\markboth{C. Ferreira and B. Parlett}{Real \textsl{dqds} for the nonsymmetric tridiagonal eigenvalue problem}

%%%%%%%%%%%%%%%%%%%%%%%%%%%%%%%%%%%%%%%%%%%%%%%%%%%%%%%%%%%%%%%%%%%%%%%%%%%%%%%%%%%
\section{Introduction}
The \textsl{dqds} algorithm was introduced in $1994$ in \cite{Fernando1} as a fast and extremely accurate way
to compute all the singular values of a bidiagonal matrix $B$. This algorithm implicitly performs the Cholesky
\textrm{LR} iteration on the tridiagonal matrix $B^TB$ and it is used in \textrm{LAPACK}. However the \textsl{dqds}
algorithm can also be regarded as executing, implicitly, the \textrm{LR} algorithm applied to any tridiagonal matrix with $1$'s
on the superdiagonal. Our interest is in real matrices which may have complex conjugate pairs of eigenvalues. It is natural to try to retain
real arithmetic and yet permit complex shifts of origin. Our analogue of the \textsl{double shift} \textrm{QR} algorithm of J. G. F. Francis
\cite{Francis1} is the \textsl{triple step dqds} algorithm. The purpose of this paper is to explain why 3 steps are needed to derive the algorithm,
to explain how we reject transforms with unacceptable element growth and to compare performance with some rival methods. Our conclusion is that this procedure is clearly the fastest method available at the present time.

We say nothing about the need for a tridiagonal eigensolver because this issue
is admirably covered in Bini, Gemignani and Tisseur \cite{Tisseur}. In fact many parts of \cite{Tisseur} have been of great help to us.
We also acknowledge the preliminary work on this problem by Z. Wu in \cite{Wu}.

We do not follow Householder conventions except that we reserve capital Roman letters for matrices. Section \ref{Section2}  describes other
methods, Section \ref{Section3} presents standard, but needed, material on \textrm{LR}, \textsl{dqds},
double shifts and the implicit L theorem. Section \ref{Section4} develops our \textsl{tridqds} algorithm, Section \ref{Section5} is our error
analysis, Section \ref{Section6} our splitting, deflation and shift strategy, and Section \ref{Section7} presents our numerical tests using
\textsc{Matlab}. Finally, Section \ref{Section8} gives our conclusions and also our ideas about why \textsl{tridqds} is only one (important)
ingredient for a procedure that must also provide condition numbers and eigenvectors.
%%%%%%%%%%%%%%%%%%%%%%%%%%%%%%%%%%%%%%%%%%%%%%%%%%%%%%%%%%%%%%%%%%%%%%%%%%%%%%%%%%%

\section{Other methods}\label{Section2}
\subsection{2 steps of \textrm{LR} = 1 step of \textrm{QR}}
\label{2LR1QR}
A frequent exercise for students is to show that for a symmetric positive definite tridiagonal matrix 2 steps of
the \textrm{LR} (Cholesky) algorithm produces the same matrix as 1 step of the \textrm{QR} algorithm. Less well known is the
article by H. Xu \cite{Xu} which extends this result when the symmetric matrix is not positive definite. The catch here is that
the \textrm{LR} transform, if it exists, does not preserve symmetry. The remedy is to regard similarities by diagonal matrices as ``trivial",
always available, operations. Indeed, diagonal similarities cannot introduce zeros into a matrix. So, when successful, 2 steps os \textrm{LR}
are diagonally similar to one step of \textrm{QR}. Even less well known is a short paper by J. Slemons \cite{Slemons}
showing that for a tridiagonal matrix, not necessarily symmetric, 2 steps of of \textrm{LR} are diagonally equivalent to 1 step of \textrm{HR},
see \cite{Angelica}. Note that when symmetry disappears then \textrm{QR} is out of the running because it does not preserve the tridiagonal property.

The point of listing these results is to emphasize that 2 steps of \textrm{LR} gives twice as many shift opportunities as 1 step of \textrm{QR} or
\textrm{HR}. Thus convergence can be more rapid with \textrm{LR} (or \textsl{dqds}) than with \textrm{QR} or \textrm{HR}. This is one of the reasons
that \textsl{dqds} is faster than \textrm{QR} for computing singular values of bidiagonals. This extra speed is an additional bonus to the fundamental advantage that \textsl{dqds} delivers high relative accuracy in all the singular values. The one drawback to \textsl{dqds}, for bidiagonals, is that
the singular values must be computed in monotone increasing order; \textrm{QR} allows the singular values to be found in any order.

In our case, failure is always possible and so there is no constraint on the order in which eigenvalues are found. The feature of having more opportunities to shift leads us to favor \textsl{dqds} over \textrm{QR} and \textrm{HR}. See the list of other methods which follows. We take up the methods in historical order and consider only those that preserve tridiagonal form.
%%%%%%%%%%%%%%%%%%%%%%%%%%%%%%%%%%%%%%%%%%%%%%%%%%%%%%%%%%%%%%%%%%%%%%%%%%%%%%%%%%%
%%%%%%%%%%%%%%%%%%%%%%%%%%%%%%%%%%%%%%%%%%%%%%%%%%%%%%%%%%%%%%%%%%%%%%%%%%%%%%%%%%%
\subsection{Cullum's complex \textrm{QR} algorithm}
As part of a program that used the Lanczos algorithm to reduce a given matrix to tridiagonal form in \cite{Cullum}, Jane Cullum used the fact
that an unsymmetric tridiagonal matrix may always be balanced by a diagonal similarity transformation. She then observed that another diagonal similarity with $1$ or $i$ produces a symmetric, but complex, tridiagonal matrix to which the (complex) tridiagonal \textrm{QR} algorithm may be applied. The process is not backward stable because the relation
\[
\cos^2\tau+\sin^2 \tau =1
\]
is not constraint on $|\cos \tau|$ and $|\sin \tau|$ when they are not real. Despite the possibility of breakdown the method proved satisfactory
in practice. We have not used it in our comparisons because we are persuaded by \ref{2LR1QR} that it is out performed by the complex \textsl{dqds} algorithm, described below.

%%%%%%%%%%%%%%%%%%%%%%%%%%%%%%%%%%%%%%%%%%%%%%%%%%%%%%%%%%%%%%%%%%%%%%%%%%%%%%%%%%%
\subsection{Liu's \textrm{HR} algorithm}
In \cite{Liu} Alex Liu found a variation on the \textrm{HR} algorithm of Angelika Bunse-Gerstner that, in exact arithmetic, is guaranteed
not to breakdown - but the price is a temporary increase in bandwith. This procedure has only been implemented in \textsc{Maple}
and we do not include it in our comparison.

%%%%%%%%%%%%%%%%%%%%%%%%%%%%%%%%%%%%%%%%%%%%%%%%%%%%%%%%%%%%%%%%%%%%%%%%%%%%%%%%%%%
\subsection{Complex \textsl{dqds}}
In his thesis David Day \cite{Day} developed a Lanczos-style algorithm to reduce a general matrix to tridiagonal
form and, as with Jane Cullum,  needed a suitable algorithm to compute its eigenvalues. He knew of the effectiveness
of \textsl{dqds} in the symmetric positive definite case and realized that \textsl{dqds} extends formally to any tridiagonal that admits
triangular factorization. Without positivity the attractive property of achieving high relative accuracy disappears but,
despite possible element growth, the error analysis for \textsl{dqds} persists: if the transform does not breakdown then tiny well chosen
changes in the entries of input $L,U$ (giving $\widetilde{L},\widetilde{U}$) and output $\hat{L},\hat{U}$ (giving $\breve{L},\breve{U}$)
produces an exact relation
\[
\breve{L}\breve{U}=\widetilde{L}\widetilde{U} -\sigma I
\]
with the given shift $\sigma$. See Section \ref{Section51}. The code uses complex arithmetic because of the possible presence of complex conjugate pairs of eigenvalues.
We have wrapped David Day's complex \textsl{dqds} code in a more sophisticated wrapper that chooses suitable shifts after rejecting a transform for
excessive element growth.

%%%%%%%%%%%%%%%%%%%%%%%%%%%%%%%%%%%%%%%%%%%%%%%%%%%%%%%%%%%%%%%%%%%%%%%%%%%%%%%%%%%
\subsection{Ehrlich-Aberth algorithm}
\label{BGT}
This very careful and accurate procedure was presented
by Bini, Gemignani and Tisseur in \cite{Tisseur}. It finds the zeros of the characteristic polynomial
$p(\cdot)$ and exploits the tridiagonal form to evaluate $p^\prime(z)/p(z)$ for any $z$.
The polynomial solver improves a full set of approximate zeros at each step. Initial approximations
are found using a divide-and-conquer procedure that delivers the eigenvalues of the top and bottom halves
of the matrix $T$. The quantity $p^\prime(z)/p(z)$ is evaluated indirectly as $\left[\trace(zI-T)^{-1}\right]$
using a \textrm{QR} factorization of $zI-T$. Since $T$ is not altered there is no deflation to
assist efficiency. Very careful tests exhibit the method's accuracy - but it is very slow compared to both
\textsl{dqds}-type algorithms.

%%%%%%%%%%%%%%%%%%%%%%%%%%%%%%%%%%%%%%%%%%%%%%%%%%%%%%%%%%%%%%%%%%%%%%%%%%%%%%%%%%%
\section{LR and \textsl{dqds}}\label{Section3}
The reader is expected to have had some exposure to the \textrm{QR} and/or \textrm{LR} algorithms so we will be brief.

%%%%%%%%%%%%%%%%%%%%%%%%%%%%%%%%%%%%%%%%%%%%%%%%%%%%%%%%%%%%%%%%%%%%%%%%%%%%%%%%%%%
\subsection{LU factorization} Any $n\times n$ matrix $A$ permits unique triangular factorization
$A=LD\widetilde{U}$ where $L$ is unit lower triangular, $D$ is diagonal, $\widetilde{U}$ is unit upper triangular,
if and only if the leading principal submatrices of orders $1,\dots,n-1$ are nonsingular.

In this paper we follow common practice and write $U=D\widetilde{U}$ so that the ``pivots" (entries of $D$)
lie on $U$'s diagonal. Throughout this paper any matrix $L$ is unit lower triangular and $U$ is upper triangular.

%%%%%%%%%%%%%%%%%%%%%%%%%%%%%%%%%%%%%%%%%%%%%%%%%%%%%%%%%%%%%%%%%%%%%%%%%%%%%%%%%%%
\subsection{LR transform with shift} Note that $U$ is ``right" triangular and $L$ is
``left" triangular and this explains the standard name \textrm{LR}. For any shift $\sigma$
let
\begin{align}
A-\sigma I&=LU,\\
\widehat{A}&=UL+\sigma I.
\end{align}
Then $\widehat{A}$ is the \textrm{LR}$(\sigma)$ transform of $A$. Note that
\[
\widehat{A}=L^{-1}(A-\sigma I)L+\sigma I = L^{-1}AL.
\]
We say that the shift is restored (in contrast to \textsl{dqds} - see below). The \textrm{LR} algorithm consists of repeated \textrm{LR} transforms with shifts chosen to enhance convergence to upper triangular form. For the theory see \cite{Rut3,Rut4,Watkins1,Watkins2}.

In contrast to the well known \textrm{QR} algorithm, the \textrm{LR} algorithm can breakdown and can suffer from element growth,
$\norm{L}>>\norm{A}$, $\norm{U}>>\norm{A}$. However \textrm{LR} preserves the banded form of $A$ while \textrm{QR} does not
(except for the Hessenberg form).

When a matrix $A$ is represented by its entries then the shift operation
$A\longrightarrow A-\sigma I$ is trivial. When a matrix is given in factored form the shift operation
is not trivial.

%%%%%%%%%%%%%%%%%%%%%%%%%%%%%%%%%%%%%%%%%%%%%%%%%%%%%%%%%%%%%%%%%%%%%%%%%%%%%%%%%%%
\subsection{The \textsl{dqds} algorithm}\label{sectiondqds}
 From now on we focus on tridiagonal matrices in $J$-form - entries $(i,i+1)$ are all $1$, $i=1,\ldots,n-1$.
Throughout this paper all $J$ matrices have this form.

If $J-\sigma I$ permits triangular factorization
\[
J-\sigma I=LU
\]
then $L$ and $U$ must have the following form
\begin{equation}
\begin{array}{cc}
  L= \begin{bmatrix}
     1   &       &              &          &    \\
           l_1   &1             &          &        &    \\
                 &\ddots        &\ddots    &        &    \\
                 &              &l_{n-2}   &1       &    \\
                 &              &          &l_{n-1} &1
  \end{bmatrix},
&\qquad
  U= \begin{bmatrix}
     u_1 &1    &       &                 &    \\
               &u_2    &1                &         &    \\
               &       &       \ddots    &\ddots   &    \\
               &       &                 &u_{n-1}  &1   \\
               &       &                 &         &u_{n}
 \end{bmatrix}.
\end{array}
\end{equation} \label{matLandU}
It is an attractive feature of \textrm{LR} that
\[
UL=\widehat{J}
\]
is also of $J$-form. Thus the parameters $l_i$, $i=1,\ldots,n-1$, and $u_j$, $j=1,\ldots,n$, determine
the matrices $L$ and $U$ above and implicitly define two tridiagonal matrices $LU$ and $UL$.

The \textsl{qds} algorithm is equivalent to the \textrm{LR} algorithm but no tridiagonal matrices are
ever formed. The \textsl{progressive}  transformation is from $L,U$ to $\widehat{L}, \widehat{U}$,
\begin{equation}
\widehat{L}\widehat{U}=UL-\sigma I \label{qdstransform}.
\end{equation}
Notice that the shift is not restored and so $\widehat{U}\widehat{L}$ is
not similar to $UL$.

Equating entries in each side of equation (\ref{qdstransform}) gives
$$
\begin{array}{rl}
\textbf{\textsl{qds}}(\sigma):& \hat{u}_1=u_1+l_1-\sigma;\\
             & \mbox{\textbf{for }}i=1,\ldots,n-1\\
             & \hspace{0.6 cm} \hat{l}_i=l_iu_{i+1}/\hat{u}_{i}\\
             & \hspace{0.6 cm} \hat{u}_{i+1}=u_{i+1}+l_{i+1}-\sigma-\hat{l}_i\\
             & \mbox{\textbf{end for}}.
\end{array}
$$

\noindent The algorithm \textsl{qds} fails when $\hat{u}_{i}=0$ for some
$i<n.$ When $\sigma=0$ we write simply \textsl{qd}, not \textsl{qds}.

In 1994 a better way was found to implement $\textsl{qds}(\sigma)$ that had been used by
Rutishauser as early as 1955. These are \textsl{called differential qd} algorithms. See
\cite{Parlett2} for more history. This form uses uses an extra variable $d$ but has compensating advantages.
\begin{equation*}
\begin{array}{rl}
\textbf{\textsl{dqds}}(\sigma):& d_1=u_1-\sigma\\
             & \mbox{\textbf{for }}i=1,\ldots,n-1 \\
             & \hspace{0.6 cm} \hat{u}_i=d_i+l_i \\
             & \hspace{0.6 cm} \hat{l}_i=l_i(u_{i+1}/\hat{u}_i)\\
             & \hspace{0.6 cm} d_{i+1}=d_i(u_{i+1}/\hat{u}_{i})-\sigma\\
             & \mbox{\textbf{end for}}\\
             & \hat{u}_n=d_n.
\end{array}
\end{equation*}
By definition, \textsl{dqd}=\textsl{dqds}(0).

A word on terminology. In Rutishauser's original work $q_i=u_i$, $e_i=l_i$; the
$q_i$'s were certain \textsl{quotients} and and the $e_i$'s were called \textsl{modified differences}.
In fact the \textsl{qd} algorithm led to the \textrm{LR} algorithm, not vice-versa. The reader can find more information concerning \textsl{dqds} in \cite{Parlett2,Parlett5}

\smallskip

One virtue of the \textsl{dqds} and \textrm{QR} transforms is that they work on the whole matrix so that
large eigenvalues are converging near the top, albeit slowly, while the small ones are being picked off at the bottom.

\smallskip

We summarize some advantages and disadvantages of the factored form.

\medskip

\noindent \textbf{Advantages of the factored form}
\smallskip
\begin{enumerate}
\item
$L,U$ determines the entries of $J$ to greater than
working-precision accuracy because the addition and
multiplication of $l$'s and $u$'s is implicit. Thus, for instance,
the $(i,i)$ entry of $J$ is given by $l_{i-1}+u_i$ implicitly but $fl(l_{i-1}+u_i)$
explicitly.
\item
Singularity of $J$ is detectable by inspection when $L$ and $U$
are given, but only by calculation from $J$. So, $LU$ reveals singularity, $J$ does not.
\item
$LU$ defines the eigenvalues better than $J$ does (usually). There is more on this in \cite{DillonParlett2}.
\item
Solution of $Jx=b$ takes half the time when $L$ and $U$ are
available.
\end{enumerate}

\bigskip

\noindent \textbf{Disadvantages of the factored form}

\smallskip
The mapping $J,\sigma \mapsto  L,U$ is not everywhere defined for all pairs $J,\sigma$ and can suffer from element
growth. This defect is not as serious as it was when the new transforms were written over the old ones.
For tridiagonals we can afford to double the storage and map
$L,U$ into different arrays $\widehat{L},\widehat{U}$. Then we can decide whether or not to accept $\widehat{L},\widehat{U}$ and only then would $L$
and $U$ be overwritten. So the difficulty of excessive
element growth has been changed from disaster to the non-trivial but less intimidating one of, after rejecting a transform,
choosing a new shift that will not spoil convergence and will not cause another rejection.
\smallskip

Now we turn to our main question of $\textsl{dqds}(\sigma)$: how can complex shifts be used without having to use complex arithmetic? This question has a beautiful answer for \textrm{QR} and \textrm{LR} iterations.

%%%%%%%%%%%%%%%%%%%%%%%%%%%%%%%%%%%%%%%%%%%%%%%%%%%%%%%%%%%%%%%%%%%%%%%%%%%%%%%%%%%
\subsection{Double shift \textrm{LR} algorithm}\label{DoubleLR}
We use the $J,L$ and $U$ notation from the previous section. Consider two steps of the
\textrm{LR} algorithm with shifts $\sigma_1$ and $\sigma_2$,

\begin{align*}
J_2-\sigma_1I  &=L_2U_2\\
          J_3  &=U_2L_2+\sigma_1I\\
J_3-\sigma_2I  &=L_3U_3\\
          J_4  &=U_3L_3+\sigma_1I.
\end{align*}

Then
\begin{equation}
J_4=\bs{\cal{L}}^{-1} J_2 \bs{\cal{L}}\label{J4}
\end{equation}
with
\[
\bs{\cal{L}}=L_2L_3,\quad  \bs{\cal{U}}=U_3U_2
\]
and
\begin{align}
\bs{\cal{L}}\bs{\cal{U}}&=L_2(J_3-\sigma_2I)U_2 \nonumber \\
 &=L_2(U_2L_2+\sigma_1I)U_2-\sigma_2L_2U_2  \nonumber\\
 &=L_2U_2\left[ L_2U_2+(\sigma_1-\sigma_2)I\right]  \nonumber\\
 &=(J_2-\sigma_1I)(J_2-\sigma_2I)  \nonumber \\
 &=J_2^2-(\sigma_1+\sigma_2)J_2+\sigma_1\sigma_2I =:M \label{matrixM}
\end{align}
Suppose that $J_2$ is real and $\sigma_1$ is complex. Then $J_4$ will be real if, and only if,
$\sigma_2=\bar{\sigma}$. The reason is that $M$ is real, so that $\bs{\cal{L}}$ and $\bs{\cal{U}}$
are real and, by (\ref{J4}), $J_4$ is the product of real matrices. Note however that $L_2,U_2,L_3,U_3$
are all complex. Fortunately it is possible to compute $J_4$ from $J_2$ without using $J_3$. This depends on the following result.

\smallskip
\begin{theorem}\label{ImplicitLtheorem}
\textsc{[Implicit L theorem]} If $H_1$ and $H_2$ are unreduced upper Hessenberg matrices and $H_2=L^{-1}H_1L$,
where $L$ is unit lower triangular, then $H_2$ and $L$ are completely determined
by $H_1$ and column 1 of $L$, $L\boldsymbol{e}_1$.
\end{theorem}
We omit the proof.

\smallskip
The clever application to $J_2$ and $J_4$ is to observe that column 1 of $M$,
\[
M\boldsymbol{e}_1=\bs{\cal{L}}\bs{\cal{U}}\boldsymbol{e}_1=\bs{\cal{L}}\boldsymbol{e}_1u_{11}, \quad u_{11}=m_{11},
\]
is proportional to column 1 of $\mathbf{L}$ and has only three nonzero entries below the diagonal because $J_2$ is tridiagonal.
Now choose
\[
\mathcal{L}_1=I+\boldsymbol{m}\boldsymbol{e_1}^T
\]
where
\[
\boldsymbol{m}=\begin{bmatrix}
               0 & m_{21}/m_{11} & m_{31}/m_{11} & 0 &\ldots &0
               \end{bmatrix}^T
\]
and perform an explicit similarity transform on $J_2$,
\[
\mathcal{L}_1^{-1}J_2\mathcal{L}_1=:K.
\] Observe that $K$ is not tridiagonal.
In the $6\times 6$ case
\[
K=\begin{bmatrix}
   x & 1 &   &   &   &  \\
   x & x & 1 &   &   &  \\
   + & x & x & 1 &   &  \\
   + &   & x & x & 1 &  \\
     &   &   & x & x & 1\\
     &   &   &   & x & x
  \end{bmatrix}.
\]
Next we apply a sequence of elementary similarity transformations such that each
transformation pushes the $2\times 1$ bulge one row down and one column to the right. Finally
the bulge is chased off the bottom to restore the $J$-form. In exact arithmetic, the implicit L theorem
ensures that this technique of \textsl{bulge chasing} gives
\[
J_4=(\mathcal{L}_1 \ldots \mathcal{L}_{n-1})^{-1}J_2(\mathcal{L}_1 \ldots \mathcal{L}_{n-1}) \;\; \text{ and } \;\;
\bs{\cal{L}}=\mathcal{L}_1 \ldots \mathcal{L}_{n-1}.
\]

%%%%%%%%%%%%%%%%%%%%%%%%%%%%%%%%%%%%%%%%%%%%%%%%%%%%%%%%%%%%%%%%%%%%%%%%%%%%%%%%%%%
\section{Triple \textsl{dqds} algorithm}\label{Section4}

\subsection{Connection to \textrm{LR} algorithm}
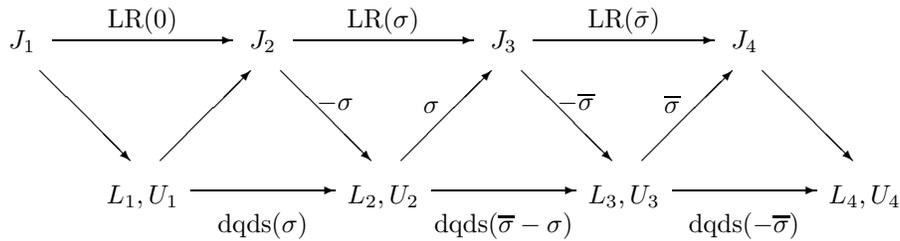
\begin{figure}[htbp] %h-here; b-Bottom; t-top; p-numa p\'{a}gina separada
\bigskip
\bigskip
\centering

\input{FiguraRelLRqdsTalkTese}
\caption{Double shift \textrm{LR} and three steps of \textsf{dqds}}
\label{figDoubleLRTridqds}

\end{figure}

In figure \ref{figDoubleLRTridqds} we examine the double shift \textrm{LR} transform derived in section \ref{DoubleLR}
but with a significant difference. Instead of $J_2$ being an arbitrary real matrix in $J$-form, we assume that it is given to us in the form
$U_1L_1$ obtained from one step of the \textrm{LR} algorithm with shift $0$ from real $J_1$.

The crucial observation is that, along the bottom line $L_2,U_2,L_3,U_3$ are all complex and so it requires 3 \textsl{dqds}
steps to go from real $L_1,U_1$ to real $L_4,U_4$. Moreover the non-restoring shifts in \textsl{dqds} are
\[
\sigma-0, \quad \bar{\sigma}-\sigma, \quad 0-\bar{\sigma}.
\]
Here is another way of seeing the relation between \text{LR} and \textsl{dqds}:
\begin{align*}
&\left\{\begin{array}{ll}
         J_1 =& L_1U_1 \\
         J_2 =& U_1L_1
        \end{array} \right.                   \\\\
&\left\{\begin{array}{rll}
          J_2-\sigma I=& L_2U_2  &  \qquad L_2U_2= U_1L_1-\sigma I \\
          J_3=& U_2L_2+\sigma I
          \end{array} \right.        \\\\
&\left\{\begin{array}{rll}
        J_3-\overline{\sigma}I=& L_3U_3  &  \qquad  L_3U_3=(U_2L_2+\sigma I)-\overline{\sigma}I \\
        J_4=& U_3L_3+\overline{\sigma} I
        \end{array} \right.            \\\\
&\left\{\begin{array}{rll}
        J_4=& L_4U_4   &  \qquad L_4U_4=(U_3L_3+\overline{\sigma}I)-0I\\
        \ldots &\phantom{J_4= U_3L_3+\overline{\sigma} I}
        \end{array} \right.
\end{align*}
Recall from the previous section that the double \textrm{LR} algorithm can work with complex shifts in real arithmetic by
\textsl{bulge chasing}. The rest of this section developes a form of bulge chasing for the \textsl{dqds} algorithm.

%%%%%%%%%%%%%%%%%%%%%%%%%%%%%%%%%%%%%%%%%%%%%%%%%%%%%%%%%%%%%%%%%%%%%%%%%
\subsection{3 steps of \textsl{dqds}} In contrast to a single \textsl{dqds}
step our \textsl{triple dqds} restores the shifts. Recall from (\ref{J4}) in section \ref{DoubleLR} that
\begin{equation}
L_4U_4=J_4=\bs{\cal{L}}^{-1} J_2 \bs{\cal{L}}=\bs{\cal{L}}^{-1} U_1L_1 \bs{\cal{L}}
\end{equation}
and, since $\sigma_1=\sigma \text{ and } \sigma_2=\bar{\sigma}$, matrix $M$ in (\ref{matrixM}) is given by
\begin{equation}
M=(U_1L_1)^2-2(\Re\sigma_1)U_1L_1+|\sigma_1|^2I.\label{M3dqds}
\end{equation}
The idea is to transform $U_1$ into $L_4$ and $L_1$ into $U_4$ by bulge chasing in each matrix,
\[
L_4U_4=\underbrace{\bs{\cal{L}}^{-1}U_1}\underbrace{L_1\bs{\cal{L}}}.
\]
Notice that we need to transform an upper bidiagonal into a lower bidiagonal and vice-versa.
From the uniqueness of the $LU$ factorization, when it exists, it follows that there is a unique hidden matrix $X$ such that
\[
L_4=\bs{\cal{L}}^{-1}U_1X^{-1}, \quad XL_1\bs{\cal{L}}=U_4.
\]
For more on $X$ see \cite{Parlett1}.
The matrix $\bs{\cal{L}}$ is given, from section \ref{DoubleLR} as a product
\[
\bs{\cal{L}}=\mathcal{L}_1\ldots\mathcal{L}_{n-1}\mathcal{L}_{n}
\]
($\mathcal{L}_{n}=I$) and we will gradually construct the matrix $X$ in corresponding factored form $X_n,\ldots,X_2X_1$.
In fact we will write each $X_i$ as a product
\[
X_i=Y_iZ_i.
\]
The details are quite complicated.

%%%%%%%%%%%%%%%%%%%%%%%%%%%%%%%%%%%%%%%%%%%%%%%%%%%%%%%%%%%%%%%%%%%%%%%%%%%%%%%%%%%
\subsubsection{Chasing the bulges}\label{chasing}
Starting with the factors $L_1$, $U_1$ and the shift $\sigma$, we
normalize column 1 of $M$ in (\ref{M3dqds}) to form ${\mathcal L}_1$,
spoil the bidiagonal form with
\[
\underbrace{\mathcal{L}_1^{-1}U_1}\underbrace{L_1\mathcal{L}_1}
\]
and at each \textsl{minor} step $i$, $i=1,\ldots,n$, matrices $Z_i$, $\mathcal{L}_i$ and $Y_i$ are chosen to
chase the bulges. After $n$ minor steps, we obtain $L_4$ and $U_4$,
\begin{align*}
L_4U_4=&\underbrace{\mathcal{L}_n^{-1}\cdots\mathcal{L}_1^{-1}U_1Z_1^{-1}Y_1^{-1}\cdots
Z_n^{-1}Y_n^{-1}} \underbrace{Y_nZ_n\cdots
Y_1Z_1L_1\mathcal{L}_1\cdots\mathcal{L}_n}\\
    =&\underbrace{\mathcal{L}_n^{-1}\cdots\mathcal{L}_1^{-1}U_1X_1^{-1}\cdots
X_n^{-1}} \underbrace{X_n\cdots
X_1L_1\mathcal{L}_1\cdots\mathcal{L}_n}\\
  =&\underbrace{\bs{\cal{L}}^{-1}U_1X^{-1}} \underbrace{XL_1\bs{\cal{L}}}
\end{align*}

Conceptually we create two work arrays $F$ and $G$. Initially,
\[
F=U_1, \qquad G=L_1
\]
and, finally,
\[
F=L_4, \qquad G=U_4.
\]
For a complex shift $\sigma$, the triple \textsl{dqds} algorithm has the following matrix formulation:

\begin{equation*}
\begin{array}{rll}
\textbf{\textsl{tridqds}}(\sigma,\bar{\sigma}):& \\
             & F=U_1;\; G=L_1\\
             & F=FZ_1^{-1};\;  G=Z_1G\\
             & F=\mathcal{L}_1^{-1}F; \; G=G\mathcal{L}_1 & [\text{form }\mathcal{L}_1 \text{ using } (\ref{M3dqds})]\\
             & F=FY_1^{-1};\; G=Y_1G \\\\
             & \textbf{for } i=2, \ldots, n-3 \\
             & \hspace{0.6 cm} F=FZ_i^{-1};\;  G=Z_iG\\
             & \hspace{0.6 cm} F=\mathcal{L}_i^{-1}F;\;G=G\mathcal{L}_i\\
             & \hspace{0.6 cm} F=FY_i^{-1};\; G=Y_iG  &[Z_i \text{ with one, } \mathcal{L}_i \text{ with two and } Y_i\text{ with three}\\
             & \textbf{end for} \hspace{3 cm}         &\text{\;nonzero off-diagonal entries}]\\\\
             & \textsf{\% \;\;\text{step } n-2} \\
             & F=FZ_{n-2}^{-1};\;  G=Z_{n-2}G\\
             & F=\mathcal{L}_{n-2}^{-1}F;\;G=G\mathcal{L}_{n-2}\\
             & F=FY_{n-2}^{-1};\; G=Y_{n-2}G  \;\;\;\; &[Y_{n-2} \text{ with two nonzero off-diagonal entries}]\\\\
             & \textsf{\%\;\;\text{step } n-1} \\
             & F=FZ_{n-1}^{-1};\;  G=Z_{n-1}G\\
             & F=\mathcal{L}_{n-1}^{-1}F;\;G=G\mathcal{L}_{n-1}\\
             & F=FY_{n-1}^{-1};\; G=Y_{n-1}G  & [Y_{n-1} \text{ and } \mathcal{L}_{n-1}\text{ with one nonzero off-diagonal entry}]\\\\
             & \textsf{\%\;\;\text{step } n} \\
             & \mathcal{L}_n=I; \; Y_n=I \\
             & F=FZ_n^{-1};\;  G=Z_nG      &  [Z_n \text{ diagonal}]\\
             & L_4=F; \;F_4=G
\end{array}
\end{equation*}
\bigskip
%%%%%%%%%%%%%%%%%%%%%%%%%%%%%%%%%%%%%%%%%%%%%%%%%%%%%%%%%%%%%%%%%%%%%%%%%%%%%%%%%%%%%%%%%%%%%%%%%%%%%
\subsubsection{Details of \textsl{tridqds}}

In this section we will go into the details of the \textsl{tridqds}
algorithm described in the previous section. Consider $L_1$ with subdiagonal entries
$\l_1,\ldots,l_{n-1}$ and $U_1$ with diagonal entries $u_1,\ldots,u_n$, as defined in
Section \ref{sectiondqds}, and consider matrices $L_4$ and $U_4$ with subdiagonal entries
$\hat{l}_1,\ldots,\hat{l}_{n-1}$ and diagonal entries $\hat{u}_1,\ldots,\hat{u}_n$, respectively.

For each iteration of \textsl{tridqds}, at the beginning of a minor step $i, \; i=2,\ldots,n-2,$ the \textsl{active}
$4\times 4$ windows of $F$ and $G$ are
\begin{align}
 F= \begin{bmatrix}
        \ddots&         &       &          &         &    \\
        \ddots&        1&      &          &         &    \\
               &\hat{l}_{i-1}   &u_i       &1        &         &    \\
               &+       &       &u_{i+1}   &1        &    \\
               &+       &       &          &u_{i+2}  &\ddots   \\
               &        &       &          &         &\ddots
  \end{bmatrix},
  \qquad
  G= \begin{bmatrix}
          \ddots & \ddots&       &        &        &    \\
                 &\hat{u}_{i-1}  &1       &        &       &    \\
                 &               &*       &        &       &    \\
                 &               &+       &1       &       &    \\
                 &               &+       &l_{i+1} &1      &    \\
                 &               &        &        &\ddots &\ddots
 \end{bmatrix}.\label{BeginStep}
\end{align}
Each minor step $i$, $i=2,\ldots,n-3$, consists of the following 3 parts.
\medskip

%%%%%%%%%%%%%%%%%%%%%%%%%%%%%%%%%%%%%%%%%%%%%%%%%%%%%%%%%%%%%%%%%%%%%%%%%%%%%%%%%%
\noindent $\underline{\textbf{Minor step \textit{i}}}$
\medskip

\begin{itemize}

\item[\textbf{a)}]
$F \longleftarrow FZ_i^{-1}$ puts 0 into $F_{i,i+1}$ and 1 into $F_{i,i}$

$G \longleftarrow Z_iG$ turns $G_{i,i+1}$ into 1

\begin{align*}
  Z_i^{-1}= \begin{bmatrix}
\ddots &           &       &       &     &\\
       &          1&       &      &       &        \\
       &          &       \frac{1}{u_i}      & -\frac{1}{u_i}               &       &       \\
       &          &               0&             1        &       &         \\
       &                     &               &       &1      &      \\
       &                     &               &       &       &\ddots \\
 \end{bmatrix},
  \qquad  \qquad
  Z_i= \begin{bmatrix}
  \ddots &      &       &       &       &    \\
         &      1      &       &       &     &    \\
         &      &      u_i     &1       &     &    \\
               &      &0       &1      &     &    \\
               &       &       &       &1    &    \\
               &       &       &       &     &\ddots
  \end{bmatrix},
\end{align*}

\begin{align*}
FZ_i^{-1}= \begin{bmatrix}
        \ddots&         &       &          &         &    \\
        \ddots&        1&      &          &         &    \\
               &\hat{l}_{i-1}   &1         &0        &         &    \\
               &+       &       &u_{i+1}   &1        &    \\
               &+       &       &          &u_{i+2}  &\ddots   \\
               &        &       &          &         &\ddots
  \end{bmatrix},
\;\;  \;\;
  Z_iG=\begin{bmatrix}
          \ddots & \ddots&       &        &        &    \\
                 &\hat{u}_{i-1}  &1       &        &       &    \\
                 &               &*       &1        &       &    \\
                 &               &+       &1       &       &    \\
                 &               &+       &l_{i+1} &1      &    \\
                 &               &        &        &\ddots &\ddots
 \end{bmatrix}.
\end{align*}

%%%%%%%%%%%%%%%%%%%%%%%%%%%%%%%%%%%%%%%%%%%%%%%%%%%%%%%%%%%%%%
\item[\textbf{b)}]

$F \longleftarrow \mathcal{L}_i^{-1}F$ puts 0 in $F_{i+1,i-1}$ and $F_{i+2,i-1}$

$G\longleftarrow G\mathcal{L}_i$ defines $\hat{u}_i$ and creates 3 nonzeros below it

\begin{align*}
\mathcal{L}_i^{-1}= \begin{bmatrix}
        \ddots&         &       &          &         &    \\
              & \ddots      &      &          &         &    \\
               &       &1       &         &         &    \\
               &       &*       &1   &        &    \\
               &       &*       &          &1  &   \\
               &        &       &          &         &\ddots
    \end{bmatrix}=I+ \boldsymbol{x}\boldsymbol{e}_i^T,
    \qquad  \qquad
    \mathcal{L}_i=I-\boldsymbol{x}\boldsymbol{e}_i^T,
\end{align*}

\begin{align*}
 \mathcal{L}_i^{-1}F= \begin{bmatrix}
        \ddots&         &       &          &         &    \\
        \ddots&        1&      &          &         &    \\
               &\hat{l}_{i-1}   &1         &        &         &    \\
               &        &*       &u_{i+1}   &1        &    \\
               &        &+      &          &u_{i+2}  &1   \\
               &        &       &          &         &\ddots
  \end{bmatrix},
  \; \;\;
  G\mathcal{L}_i= \begin{bmatrix}
          \ddots & \ddots&       &        &        &    \\
                 &\hat{u}_i       &1        &       &    \\
                 &+       &1       &       &    \\
                 &+       &l_{i+1} &1      &    \\
                 &+        &        &l_{i+2} &1\\
                 &               &        &        &\ddots &\ddots
 \end{bmatrix}.
\end{align*}

\smallskip

%%%%%%%%%%%%%%%%%%%%%%%%%%%%%%%%%%%%%%%%%%%%%%%%%%%%%%%%%%%%
\item[\textbf{c)}]

$G\longleftarrow Y_iG$ puts 0 in $G_{i+1,i}$, $G_{i+2,i}$ and $G_{i+3,i}$

$F \longleftarrow FY_i^{-1}$ creates $\hat{l}_i$ and puts 2 nonzeros below it

\begin{align*}
   Y_i^{-1}= \begin{bmatrix}
          \ddots &           &   &        &        &    \\
                 &1          &        &       &    \\
                 &*        &1       &       &    \\
                 &*        &  &1      &    \\
                 &*        &        &  &1\\
                 &               &        &        & &\ddots
  \end{bmatrix}=I+ \boldsymbol{y}\boldsymbol{e}_i^T,
    \qquad  \qquad
    Y_i=I-\boldsymbol{y}\boldsymbol{e}_i^T,
  \end{align*}

\begin{align*}
 FY_i^{-1}= \begin{bmatrix}
        \ddots&         &       &          &         &    \\
        \ddots&        1&      &          &         &    \\
               &\hat{l}_i   &u_{i+1}       &1        &         &    \\
               &+       &       &u_{i+2}   &1        &    \\
               &+       &       &          &u_{i+3}  &\ddots   \\
               &        &       &          &         &\ddots
  \end{bmatrix},
  \; \;
  Y_iG= \begin{bmatrix}
          \ddots & \ddots&       &        &        &    \\
                 &\hat{u}_i  &1       &        &       &    \\
                 &               &*       &        &       &    \\
                 &               &+       &1       &       &    \\
                 &               &+       &l_{i+2} &1      &    \\
                 &               &        &        &\ddots &\ddots
 \end{bmatrix}.\\
\end{align*}

%%%%%%%%%%%%%%%%%%%%%%%%%%%%%%%%%%%%%%%%%%%%%%%%%%%%%%%%%%%%
%%%%%%%%%%%%%%%%%%%%%%%%%%%%%%%%%%%%%%%%%%%%%%%%%%%%%%%%%%%%
\end{itemize}

The result of this minor step is that the active windows of $F$ and $G$ shown in
(\ref{BeginStep}) have been moved down and to the right by one place.
See Appendix \ref{ApendiceA} for more details on the practical implementation.

Naturally steps $1, n-2, n-1,n$ are slightly different and may be found on pp. 147-157 of \cite{CFerrThesis}.
%%%%%%%%%%%%%%%%%%%%%%%%%%%%%%%%%%%%%%%%%%%%%%%%%%%%%%%%%%%%
\subsection{Operation count for \textsl{tridqds}}
In this section we will see how three steps of simple \textsl{dqds} algorithm compares
with one step of \textsl{tridqds} in what respects to the number of floating point operations required.

Here is the inner loop of \textsl{tridqds}. See Appendix \ref{ApendiceB}.
$$
\begin{array}{rll}
\textbf{\textsl{tridqds}}(\sigma,\bar{\sigma}):& \\
             & \textbf{for } i=2, \ldots, n-3 \\   %FOR
             & \hspace{0.6 cm} x_r =x_r*u_i+y_r   \\
             & \hspace{0.6 cm} x_l= -x_l*(1/\hat{l}_{i-1}) ; \;\; y_l=-y_l*(1/\hat{l}_{i-1});\;\;\\
             & \hspace{0.6 cm} \hat{u}_i=x_r-x_l;\\
             & \hspace{0.6 cm} x_r=y_r-x_l; \;\; y_r=z_r-y_l-x_l*l_{i+1}; \;\; z_r=-y_l*l_{i+2}\\
             & \hspace{0.6 cm} x_r = x_r*(1/\hat{u}_i); \;\; y_r = y_r*(1/\hat{u}_i); \;\; z_r = z_r*(1/\hat{u}_i)\\
             & \hspace{0.6 cm} \hat{l}_i = x_l+y_r+x_r*u_{i+1}\\
             & \hspace{0.6 cm} x_l = y_l+z_r+ y_r*u_{i+2}; \;\; y_l = z_r*u_{i+3}\\
             & \hspace{0.6 cm} x_r = 1-x_r; \;\; y_r = l_{i+1}-y_r; \;\; z_r =-z_r\\
             & \textbf{end for}\\    %END FOR

\end{array}
$$
A good compiler recognizes common subexpressions.

In contrast,
\begin{equation*}
\begin{array}{rl}
\textbf{\textsl{dqds}}(\sigma):& d_1=u_1-\sigma\\
             & \mbox{\textbf{for }}i=1,\ldots,n-1 \\
             & \hspace{0.6 cm} \hat{u}_i=d_i+l_i \\
             & \hspace{0.6 cm} \hat{l}_i=l_i(u_{i+1}/\hat{u}_i)\\
             & \hspace{0.6 cm} d_{i+1}=d_i(u_{i+1}/\hat{u}_{i})-\sigma\\
             & \mbox{\textbf{end for}}\\
             & \hat{u}_n=d_n.
\end{array}
\end{equation*}
In practice, each $d_{i+1}$ may be written over its
predecessor in a single variable $d$ and, if the common
subexpression $u_{i+1}/\hat{u}_i$ is recognized, then only one
division is needed if we use an auxiliary variable.

\smallskip
Table \ref{tripOper} below shows that the operation count of one
step of \textsl{tridqds} is comparable to three
steps of \textsl{dqds} (table expresses only the number of
floating point operations in the inner loops).

\medskip
\begin{table}[htbp]
\centering
\begin{tabular}{c||c|c|}
%\multicolumn{1}{c}{}     & \multicolumn{2}{c}{\textsf{eigval}}       & \multicolumn{2}{c}{\textsf{eig}}
\multicolumn{1}{c}{} & \multicolumn{1}{c}{\textsl{tridqds}} & \multicolumn{1}{c}{3 \textsl{dqds} steps}\\\hline \hline

   Divisions & 2 & 3 \\

   Multiplications & 11 & 6 \\

   Additions & 5 & 3 \\

   Subtractions & 6 & 3 \\

   Assignments & 16 & 12 \\

   Auxiliary variables & 5 & 2 \\ \hline

\end{tabular}
\bigskip
\caption{Operation count of \textsl{tridqds} and 3 \textsl{dqds} steps}

\label{tripOper}

\end{table}
But to make three steps of \textsl{dqds} equivalent to \textsl{tridqds} we have to consider \textsl{dqds}
in complex arithmetic and the total cost is raised by a factor of about 4. Thus, in complex arithmetic,
three steps of \textsl{dqds} are much more expensive than one step of \textsl{tridqds}.

%%%%%%%%%%%%%%%%%%%%%%%%%%%%%%%%%%%%%%%%%%%%%%%%%%%%%%%%%%%%%%%%%%%%%%%%%%%%%%%%%%%
\section{Error analysis}\label{Section5}
We turn to the effect of finite precision arithmetic on our algorithms. First consider the
\text{dqds} algorithm.
%%%%%%%%%%%%%%%%%%%%%%%%%%%%%%%%%%%%%%%%%%%%%%%%%%%%%%%%%%%%%%%%%%%%%%%%%%%%%%%%%%%
\subsection{\textsl{dqds}}
\label{Section51}
In the absence of over/underflow the algorithm enjoys the so-called \textsl{mixed relative stability} property.
\begin{theorem}
Let $\textsl{dqds}(\sigma)$ map $L,U$ into computed $\widehat{L},\widehat{U}$ with no division by zero, over/underflow.
Then well chosen small relative changes in the entries of both input and output matrices, of at most 3 ulps each,
produces new matrices, one pair mapped into the other, in exact arithmetic, by $\textsl{dqds}(\sigma)$.
\end{theorem}

See the diagram in Figure \ref{ComutDiag}. The remarkable feature here is that huge element growth does not impair the result.
However this useful property does not guarantee that \textsl{dqds} returns accurate eigenvalues.
See \cite{Fernando1,Parlett2}. For that, an extra requirement is needed such as positivity of all the parameters $u_j$, $l_j$
in the computation. This is the case for the eigenvalues of $B^TB$ where $B$ is upper bidiagonal.

What can be said in our case? We quote a result that is established by Yao Yang in his dissertation \cite{YaoYang}
and appears in \cite{Parlett2}. The clever idea is not to look at the $L$ and $U$ separately but to study their exact product
$J=LU$.

\begin{figure}[ht]

\centering

\input{CommutativeDiagrame}
$\breve{L}\breve{U}=\widetilde{U}\widetilde{L}-\sigma I$
\medskip

\caption{Effects of roundoff for \textmd{dqds}}
\label{ComutDiag}
\end{figure}

\smallskip
\begin{theorem}{[Y. Yang]} If $\textsl{dqds}(\sigma)$ maps $L,U$ into $\widehat{L},\widehat{U}$
 (with no division by 0, overflow/underflow) in the standard model of floating point arithmetic
 then there is a unique pair $\mathring{L}, \mathring{U}$ such that, in exact arithmetic, $\textsl{dqds}(\sigma)$
 maps $\mathring{L}, \mathring{U}$ into $\widehat{L},\widehat{U}$. Moreover, the associated tridiagonals satisfy,
 element by element,
\begin{gather*}
|\offdiag(\mathring{J})-\offdiag(J)|<2 \varepsilon |\offdiag(J)|\\
|\diag(\mathring{J})-\diag(J)|<
\varepsilon\left(2|\boldsymbol{u}|+|\sigma||\boldsymbol{1}|+|\boldsymbol{\hat{l}}|+|\boldsymbol{\hat{u}}|+2|\boldsymbol{d}|\right)
\end{gather*}
\smallskip
where $\varepsilon$ is the \textsl{roundoff unit}.
\end{theorem}

This result is Corollary 3 in Section 9 of \cite{Parlett2}. It shows that it is only the diagonal of $J$ that suffers large backward
error in the case of element growth. Since $\boldsymbol{\hat{u}}=\boldsymbol{d}+\boldsymbol{\hat{l}}$ the last inequality may be written as
\[
|\diag(\mathring{J})-\diag(J)|<
\varepsilon\left(2|\boldsymbol{u}|+|\boldsymbol{l}|+|\sigma||\boldsymbol{1}|+|\boldsymbol{\hat{l}}|+3|\boldsymbol{d}|\right).
\]
Recall that $d_i^{-1}=\left[(UL)^{-1}\right]_{ii},\; i=1,\ldots,n$. Thus the indices vulnerable to large backward error belong to any very small entries $\left[(UL)^{-1}\right]_{ii}$. For this reason we reject $\widehat{L},\widehat{U}$ when, element by element,
\begin{equation}
|\sigma||\boldsymbol{1}|+|\boldsymbol{\hat{l}}|+3|\boldsymbol{d}|>1000(|\boldsymbol{u}|+|\boldsymbol{l}|).\label{Reject}
\end{equation}
Recall that the error analysis is worst case. Recall also that the effect of a tiny $\hat{u}_k$ disappears for $i>k+1$.

%%%%%%%%%%%%%%%%%%%%%%%%%%%%%%%%%%%%%%%%%%%%%%%%%%%%%%%%%%%%%%%%%%%%%%%%%%%%%%%%%%%
\subsection{\textsl{tridqds}}
There are too many intermediate variables in this algorithm to permit a successful mixed error analysis.
However each minor step in the algorithm consists of 3 elementary similarity transformations on work matrices
$F,G$ or $G,F$. See ... in Section \ref{chasing}. Recall that an elementary matrix here is of the form
$I+\boldsymbol{v}\boldsymbol{e}_j^T$, with inverse $I-\boldsymbol{v}\boldsymbol{e}_j^T$, and $\boldsymbol{v}$
has at most 3 nonzero entries. So we examine the condition number of these 3 similarity transforms. Consult
Appendix \ref{ApendiceA} to follow the details.

\begin{itemize}
\item The active part of $Z_i$ is
\[
\begin{bmatrix}
u_i & 1\\
0 & 1
\end{bmatrix} \qquad \text{and} \qquad \cond(Z_i)\simeq \max\left\{|u_i|, |u_i|^{-1}\right\}.
\]
%%%%%%%%%%%%%%%%%
\item The active part of $\mathcal{L}_i$ is
\[
\begin{bmatrix}
1 & \\
-x_l/\hat{l}_{i-1} & 1\\
-y_l/\hat{l}_{i-1} & 0 & 1
\end{bmatrix}  \qquad \text{and} \qquad
\cond(\mathcal{L}_i)\simeq 1+ \left(\frac{x_l}{\hat{l}_{i-1}}\right)^2 + \left(\frac{y_l}{\hat{l}_{i-1}}\right)^2.
\]
%%%%%%%%%%%%%%%%
\item The active part of $Y_i$ is
\[
\begin{bmatrix}
1 & \\
-x_r/\hat{u}_{i} & 1\\
-y_r/\hat{u}_{i} & 0 & 1\\
-z_r/\hat{u}_{i} & 0 & 0 &1\\
\end{bmatrix}  \qquad \text{and} \qquad
\cond({Y}_i)\simeq 1+ \left(\frac{x_r}{\hat{u}_{i}}\right)^2 +
\left(\frac{y_r}{\hat{u}_{i}}\right)^2+\left(\frac{z_r}{\hat{u}_{i}}\right)^2.
\]
\end{itemize}
The variables $x_l,y_l,x_r,y_r,z_r$ are formed from additions and multiplications of previous quantities.
Note that $u_i$ is part of the input and so is assumed to be of acceptable size. We see that it is tiny values
of $\hat{l}_{i-1}$ and $\hat{u}_i$ that lead to an ill-conditioned similarity at minor step $i$.
In the simple \textsl{dqds} algorithm a small value of $\hat{u}_i$ (relative to $u_{i+1}$) leads to a large value
of $\hat{l}_{i+1}$ and $d_{i+1}$. In \textsl{tridqds} the effect of 3 consecutive transforms is more complicated.
The message is the same: reject any transform that has more then modest element growth, as determined by (\ref{Reject})
in the previous section. This challenge calls for further study.

%%%%%%%%%%%%%%%%%%%%%%%%%%%%%%%%%%%%%%%%%%%%%%%%%%%%%%%%%%%%%%%%%%%%%%%%%%%%%%%%%%%
\section{Implementation details}\label{Section6}
%%%%%%%%%%%%%%%%%%%%%%%%%%%%%%%%%%%%%%%%%%%%%%%%%%%%%%%%%%%%%%%

\subsection{Deflation $\boldsymbol{(n \leftarrow n-1)}$}
Some of out criteria for deflating come from \cite{Parlett5}, others are new.
Consider both matrices $UL$ and $LU$  and the trailing $2\times 2$ blocks,
\[
\begin{bmatrix}
l_{n-1}+u_{n-1} & 1\\
l_{n-1}u_{n}  & u_n
\end{bmatrix}, \qquad \qquad
\begin{bmatrix}
l_{n-2}+u_{n-1} & 1\\
l_{n-1}u_{n-1}  & l_{n-1}+u_n
\end{bmatrix}.
\]
Deflation $(n \leftarrow n-1)$ removes $l_{n-1}$ as well as $u_n$. Looking at entry $(n-1,n-1)$ of $UL$ shows that
a necessary condition is that $l_{n-1}$ be negligible compared to $u_{n-1}$,
\begin{equation}
|l_{n-1}|<tol\cdot|u_{n-1}|, \label{Deflatel1}
\end{equation}
for a certain tolerance $tol$ close to roundoff unit $\varepsilon$.

The $(n,n)$ entries of $UL$
and $LU$ suggest either $u_n+Acshift$ or $l_{n-1}+u_n+Acshift$ as eigenvalues.
$Acshift$ is the accumulated shift (recall that \textsl{dqds} is a non-restoring transform).
To make these consistent we require that
\begin{equation}
|l_{n-1}|<tol\cdot|u_n+Acshift|.\label{Deflatel2}
\end{equation}

Finally we must consider the change $\delta \lambda$ in the eigenvalue $\lambda$ caused by setting $l_{n-1}=0$.
We estimate $\delta \lambda$ by starting from $UL$  with $l_{n-1}=0$ and then allowing $l_{n-1}$ to grow.
To this end let $J$ be $UL$  with $l_{n-1}=0$ and $(u_n,\boldsymbol{y}^T,\boldsymbol{x})$ be the eigentriple for $J$.
Clearly $\boldsymbol{y}^T=\boldsymbol{e}_n^T$.
Now we consider perturbation theory with parameter $l_{n-1}$. The perturbing matrix $\delta J$, as $l_{n-1}$ grows, is
\[
l_{n-1}(\boldsymbol{e}_{n-1}+\boldsymbol{e}_{n}u_n)\boldsymbol{e}_{n-1}^T.
\]
By first order perturbation analysis
\[
|\delta \lambda|=\frac{|\boldsymbol{y}^T\delta J\boldsymbol{x}|}{\norm{\boldsymbol{x}}_2\norm{\boldsymbol{y}}_2}
\]
and $\norm{\boldsymbol{y}}_2=1$ in our case. So,
\[
|\delta \lambda|=
\frac{\left|l_{n-1}\boldsymbol{e}_n^T(\boldsymbol{e}_{n-1}+\boldsymbol{e}_{n}u_n)\boldsymbol{e}_{n-1}^T\boldsymbol{x}\right|}{\norm{\boldsymbol{x}}_2}
=\frac{|l_{n-1}u_n||x_{n-1}|}{\norm{\boldsymbol{x}}_2}
\]
and we use the crude bound
$
\ds\frac{|x_{n-1}|}{\norm{\boldsymbol{x}}_2}<1.
$
So, we let $l_{n-1}$ grow until the change
\begin{equation*}
|\delta \lambda|<|l_{n-1}u_n|
\end{equation*}
in eigenvalue $\lambda=u_n$ is no longer acceptable. Our condition for deflation is then
\begin{equation}
|l_{n-1}u_n|<tol \cdot |Acshift+u_n|.\label{Deflatel3}
\end{equation}

A similar first order perturbation analysis for $LU$ with $l_{n-1}=0$ will give our last condition for deflation.
For the eigentriple $(u_n,\boldsymbol{y}^T,\boldsymbol{x})$ we also have $\boldsymbol{y}^T=\boldsymbol{e}_n^T$.
The perturbing matrix is now
\[
l_{n-1}\boldsymbol{e}_{n}\left(\boldsymbol{e}_{n-1}^Tu_{n-1}+\boldsymbol{e}_{n}^T\right)
\]
and
\[
|\delta \lambda|=
\frac{\left|l_{n-1}\boldsymbol{e}_n^T\boldsymbol{e}_{n}(\boldsymbol{e}_{n-1}^Tu_{n-1}+\boldsymbol{e}_{n}^T)\boldsymbol{x}\right|}{\norm{\boldsymbol{x}}_2}
=|l_{n-1}|\frac{|u_{n-1}x_{n-1}+x_n|}{\norm{\boldsymbol{x}}_2}
<|l_{n-1}|\left(|u_{n-1}|+1\right).
\]
Finally we require
\begin{equation}
|l_{n-1}|\left(|u_{n-1}|+1\right)<tol \cdot |Acshift+u_n|.\label{Deflatel4}
\end{equation}

%%%%%%%%%%%%%%%%%%%%%%%%%%%%%%%%%%%%%%%%%%%%%%%%%%%%%%%%%%%%%%%
\subsection{Splitting and deflation $\boldsymbol{(n \leftarrow n-2)}$}
Recall that the implicit L theorem was invoked to justify the \textsl{tridqds} algorithm.
This result fails if any $l_k$, $k<n-1$ vanishes.
Consequently, checking for negligible values among the $l_k$ is a necessity, not a luxury for
increased efficiency. Consider $J=UL$ in block form
\[
  \begin{pmat}[{...|..}]
           &       &    &          \cr
           && J_1   &    &          \cr
           &       &    &&1         \cr\-
           &       &&\mu &          \cr
           &       &    &   &&J_2 & & \cr
           &       &    &          \cr
  \end{pmat}
\]
where $\mu=u_{k+1}l_k$, $k<n-1$. We can replace $\mu$ by $0$ when
\[
\spectrum(J_1) \cup \spectrum(J_2)=\spectrum(J),  \quad \textsl{to working accuracy}.
\]
However we are not going to estimate the eigenvalues of $J_1$ and $J_2$. Instead we create a local
criterion for splitting at $(k+1,k)$ as follows. Focus on the principal $4\times 4$ window of $J$ given by
\[
\begin{pmat}[{.|..}]
u_{k-1}+l_{k-1}    & 1\cr
u_kl_{k-1}         & u_{k}+l_{k}    & 1\cr\-
&  u_{k+1}l_k      & u_{k+1}+l_{k+1}    & 1\cr
& & u_{k+2}l_{k+1} & u_{k+2}+l_{k+2}\cr
\end{pmat}.
\]
Now $J_1$ and $J_2$ are both $2\times 2$ and our local criterion is
\begin{equation}
\det(J_1) \cdot \det(J_2)=\det(J), \quad \textsl{to working accuracy}.
\end{equation}
Let us see what this yields. Perform block factorization on $J$ and note
that the Schur complement of $J_1$ in $J$ is
\[
J_2^\prime=J_2-
\begin{bmatrix}
0 &\mu\\
0 &0
\end{bmatrix}
J_1^{-1}\begin{bmatrix}
0 &0\\
1 &0
\end{bmatrix}
\]
with
\begin{align*}
J_1^{-1}&=\frac{1}{det_1}\begin{bmatrix}
                 u_{k}+l_{k} & -1\\
                -u_kl_{k-1} &u_{k-1}+l_{k-1}
          \end{bmatrix}
\end{align*}
where
\[
det_1=\det(J_1)=u_{k-1}(u_k+l_k)+l_{k-1}l_k.
\]
Thus
\[
J_2^\prime=\begin{bmatrix}
u_{k+1}l_k      & u_{k+1}+l_{k+1}\\
u_{k+2}l_{k+1} & u_{k+2}+l_{k+2}
\end{bmatrix}
-\begin{bmatrix}
                \mu(u_{k-1}+l_{k-1})/det_1 &0\\
                0                          & 0
                \end{bmatrix}.
\]
Since $\det$ is linear by rows
\[
\det(J_2)-\det(J_2^\prime)=\mu(u_{k-1}+l_{k-1})(u_{k+2}+l_{k+2})/det_1.
\]
Our criterion reduces to splitting only when
\[
\det(J_2^\prime) =\det(J_2),  \quad \text{to working accuracy}.
\]
Thus we require
\begin{equation*}
|l_ku_{k+1}(u_{k+2}+l_{k+2})(u_{k-1}+l_{k-1})/det_1|<tol \cdot \left|\det(J_2)\right|.
\end{equation*}
Since
\[
det_2=\det(J_2)=u_{k+1}(u_{k+2}+l_{k+2})+l_{k+1}l_{k+2},
\]
the criterion for splitting $J$ at $(k+1,k)$ is then
\begin{equation}
|l_ku_{k+1}(u_{k+2}+l_{k+2})(u_{k-1}+l_{k-1})|<tol \cdot \left|det_1det_2\right|. \label{Split1}
\end{equation}
Finally, to remove $l_{k}$ we also need $l_k$ to be negligible compared to $u_{k}$,
\begin{equation}
|l_{k}|<tol\cdot|u_{k}|. \label{Split2}
\end{equation}

\medskip
\noindent\textbf{Deflation} $\boldsymbol{(n \leftarrow n-2)}$
\medskip

We use the same criterion for deflation $(n \leftarrow n-2)$,
but because $l_{k+2}=l_n=0$ there is a common factor $det_2$ on each side of (\ref{Split1}).
Deflate the trailing $2\times 2$ submatrix when
\begin{equation}
|l_{n-2}|<tol \cdot |u_{n-2}| \label{Deflate221}
\end{equation}
and
\begin{equation}
|l_{n-2}(u_{n-3}+l_{n-3})|<tol \cdot \left|u_{n-3}(u_{n-2}+l_{n-2})+l_{n-3}l_{n-2}.\right|. \label{Deflate222}
\end{equation}

We omit the role of $Acshift$ here because it makes the situation more complicated. We have to recall that
\textsl{tridqds} uses restoring shifts and $Acshift$ is always real. So, for complex shifts, $det_2$ is not going to zero.
In fact
\[
|det_2| \geq |\Im(\lambda)|^2
\]
where $\lambda$ is an eigenvalue of $J_2$.

When $n=3$ these criteria simplify a lot. Both reduce to
\[
|l_1|<tol \cdot |u_1|.
\]
%%%%%%%%%%%%%%%%%%%%%%%%%%%%%%%%%%%%%%%%%%%%%%%%%%%%%%%%%%%%%%%%%%%%%%%%%%%%%%%%%%%
\subsection{Shift strategy}
\label{shiftStrategy}
Although \textsl{tridqds} may be, and has been, used to compute all the eigenvalues, it seems sensible to include real
$\textsl{dqds}(\sigma)$ so that when all eigenvalues are real \textsl{tridqds} need not be called.

As with \textrm{LR}, the \textsl{dqds} algorithm with no shift gradually forces large entries to the top and brings
small entries towards the bottom. Before every transform both $l_{n-1}$ and $l_{n-2}$ are inspected.
If
\[
|l_{n-1}|<\frac{1}{2^4} \qquad \text{and} \qquad |l_{n-2}|<\frac{1}{2^4}
\]
then the code executes \textsl{dqds} transform with the \textsl{Wilkinson} shift or a \textsl{3dqds} transform with
\textsl{Francis} shifts depending on the sign of the discriminant.

An unexpected reward for having both transforms available is to cope with a rejected transform. Our strategy is simply to use the
other transform with the current shift. More precisely, given a complex shift $\sigma$, if  $\textsl{tridqds}(\sigma,\bar{\sigma})$
is rejected we try  $\textsl{dqds}(u_n)$; if for real $\tau$,  $\textsl{dqds}(\tau)$  is rejected, we try
$\textsl{tridqds}(\tau,\bar{\tau})$. So far, this has not failed.

More generally, an increase in the imaginary part of the shift increases diagonal dominance. At the extreme, consider a pair of pure
imaginary shifts $i\mu,-i\mu$, $\mu$ positive. The \textsl{tridqds} wants $UL+\mu^2I$ to permit triangular factorization. The bigger is $\mu$
the better.

\smallskip
A great attraction of \textsc{IEEE} arithmetic standard is that it allows the symbols \textsf{inf} and \textsf{NaN}.
Thus there is no need for time consuming with if statements in the main loop. At the end of the loops we test for rejection
or excessive element growth. We record the number of rejections.

%%%%%%%%%%%%%%%%%%%%%%%%%%%%%%%%%%%%%%%%%%%%%%%%%%%%%%%%%%%%%%%%%%%%%%%%%%%%%%%%%%%
\section{Numerical Examples}\label{Section7}
%%%%%%%%%%%%%%%%%%%%%%%%%%%%%%%%%%%%%%%%%%%%%%%%%%%%%%%%%%%%%%%%%%%%%%%%%%%%%%%%%%%

Those who work with well defined problems have the habit of determining the ``true"
(or most accurate) solutions and comparing computed values with them to give the error. The
condition number of every eigenvalue of a real symmetric matrix is 1, but only in the absolute sense.
The relative condition number can vary. In our case even the absolute condition numbers can rise to
$\infty$ and little is known about relative errors.

We have discovered \cite{Cferreira2} that more often than not the eigenvalues of tridiagonal, and reasonable
well balanced, matrices are well determined by an $LU$ or $\Delta LDL^T$ representation
($\Delta$ is a signature matrix). This is good news but much work remains. Our main focus is on the time it takes
to get reasonable approximations, recognizing that we do not know how well the data defines the eigenvalues.

\smallskip

We refer to the Ehrlich-Aberth algorithm (see Section \ref{BGT}) as $BGT$ and to
our code simply as \textsl{tridqds}, although we combine \textsl{tridqds} with
real \textsl{dqds} as described in Section \ref{shiftStrategy}.

Since we compare \textsc{Matlab} versions of all the codes we acknowledge that the elapsed times
are accurate to only about 0.02 seconds. However this is good enough to show the ratios between
$BGT$ and the \textsl{dqds} codes. The efficiency of complex \textsl{dqds} is harder to determine.
Sometimes the same, sometimes \textsl{tridqds} is $3$ times faster.

Since the number of iterations needed for convergence on our (modest) test bed has remained about $4n$, we have not
tried for a strategy as sophisticated as the one in \cite{Parlett5}.

%The experiments were carried out with \textsc{Matlab} 7.5.0 (R2007b) on a PC M 430 2.27GHz
%with IEEE double precision arithmetic under Windows 7.

\bigskip

\noindent{\textbf{Bessel matrix}}

\bigskip

Bessel matrices, associated with generalized Bessel polynomials, are nonsymmetric tridiagonals
matrices defined by
$B_n^{(a,b)}=\textmd{tridiag}(\boldsymbol{\beta},\boldsymbol{\alpha},\boldsymbol{\gamma})$
with
$$
\alpha_1=-\frac{b}{a},\;\;\;\gamma_1=-\alpha_1,\;\;\;
\beta_1=\frac{\alpha_1}{a+1},
$$
and
\begin{align*}
\alpha_j &:=-b\frac{a-2}{(2j+a-2)(2j+a-4)},\;\;\;j=2,\ldots,n,\\
\gamma_j &:=b\frac{j+a-2}{(2j+a-2)(2j+a-3)}, \\
\beta_j &:=-b\frac{j}{(2j+a-1)(2j+a-2)},\;\;\;j=2,\ldots,n-1.
\end{align*}
Parameter $b$ is a scaling factor and most authors take $b=2$ and so do we. The case $a\in \matR$ is the most investigated in literature.
The eigenvalues of $B_n^{(a,b)}$, well separated complex eigenvalues, suffer from ill-conditioning that increases with $n$
- close to a defective matrix. In Pasquini \cite{Pasquini} it is mentioned that the ill-conditining seems to reach its maximum when $a$ ranges from $-8.5$ to $-4.5$.

\begin{figure*}[htbp]
   \centering{
   \subfigure[$n=18$]{
   \includegraphics[width=7.6cm]{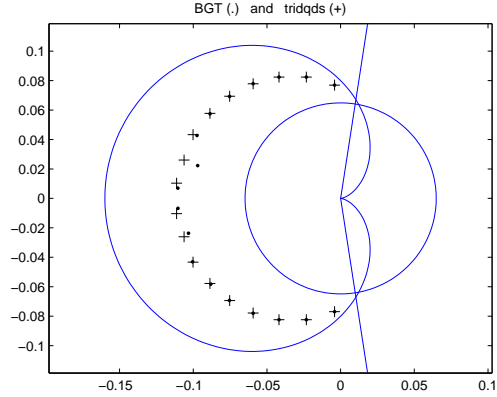}
   \label{FigBessel18}
    }}
   \centerline{
   \subfigure[$n=20$]{\includegraphics[width=7.6cm]{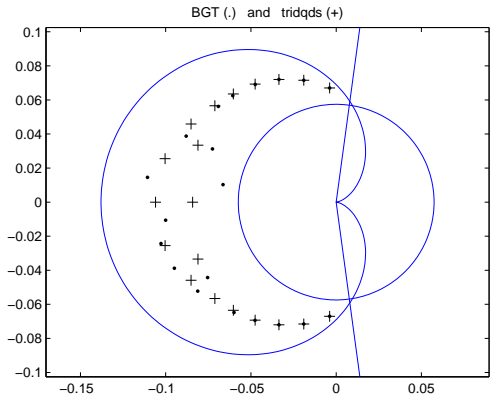}
   \label{fig:XXX}}
   \hfil
   \subfigure[$n=40$]{\includegraphics[width=7.6cm]{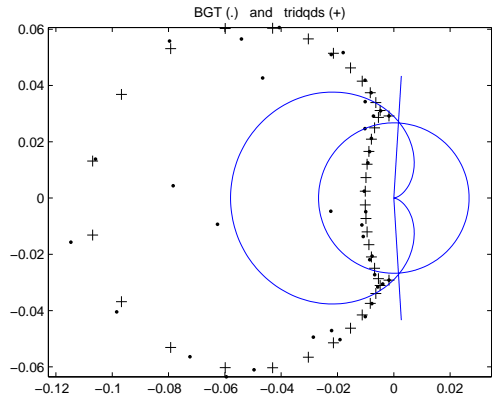}
   \label{fig:XXXX}}
   }
   \caption{\it Eigenvalues of Bessel matrix $B_{n}^{(-4.5,2)}$}
   \label{FigureBessel}
\end{figure*}

%%%%%%%%%%%%%%%%%%%%%%%%%%%%%%%%%%%%%%%%%%%%%%%%%%%%%%%%%%%%%%%%
\smallskip
Our examples take $B_{n}^{(-4.5,2)}$ for $n=18,20,40$. We show pictures for $BGT$ and
\textsl{tridqds} to illustrate the extreme sensitivity of some of the eigenvalues. The results of complex \textsl{dqds}
are visually identical to \textsl{tridqds}, so we don't show them. In exact arithmetic the spectrum lies on an arc in the interior
of the moon-shaped region. See Figure \ref{FigureBessel}.

%%%%%%%%%%%%%%%%%%%%%%%%%%%%%%%%%%%%%%%%%%%%%%%%%%%%%%%%%%%%%%%%%%%%%%%%%%%%%%%%%%%
%%%%%%%%%%%%%%%%%%%%%%%%%%%%%%%%%%%%%%%%%%%%%%%%%%%%%%%%%%%%%%%%%%%%%%%%%%%%%%%%%%%
\pagebreak

\bigskip

\noindent\textbf{Clement matrix}

\bigskip

The so-called \textsl{Clement} matrices (see \cite{Clement})
$$
C=\tridiag(\bs{b},\bs{0},\bs{c})
$$
with $b_j=j$ and $c_j=b_{n-j}$, $j=1,\ldots,n-1$,
have real eigenvalues
\begin{align*}
\pm \; n-1, n-3, \ldots, 1, \qquad \text{for } n \text{ even},\\
\pm \; n-1, n-3, \ldots, 0, \qquad \text {for } n \text{ odd}.
\end{align*}
These matrices posed no serious difficulties. The initial zero diagonal obliges the
\textsl{dqds} based methods to take care when finding an initial $L,U$ factorization.

The \textsl{tridqds} code uses only real \textsl{dqds} transforms as it should.
Our accuracy is less than $BGT$ but satisfactory.
The complex \textsl{dqds} and \textsl{tridqds} performed identically.
The ratio of elapsed times is the striking feature.

Our numerical tests have $n=100,200,400,800$. The minimum and maximum relative errors, $rel_{min}$ and $rel_{max}$,
are shown in Table \ref{tabelaRelClem} and the CPU times in Table \ref{ClementAcc}.

\begin{table}[htbp]
\centering
\begin{tabular}{l||c|c|c|c|c|c|}
\multicolumn{1}{c}{} & \multicolumn{2}{c}{$BGT$}         & \multicolumn{2}{c}{\textsl{complex dqds}} & \multicolumn{2}{c}{\textsl{tridqds}}  \\\hline
$n$ & $rel_{min}$       & $rel_{max}$       & $rel_{min}$          & $rel_{max}$     & $rel_{mim}$     & $rel_{max}$     \\ \hline \hline
100 & $0$    & $3\;10^{-16}$    & $0$       & $3\;10^{-14}$   & $0$   & $3\;10^{-14}$   \\
200 &$0$     &$4 \;10^{-16}$   & $0$        & $3\;10^{-13}$   & $0$   & $3\;10^{-13}$   \\
400 &$0$     &$1 \;10^{-15}$   & $0$        & $3\;10^{-12}$   & $0$   & $3\;10^{-12}$  \\
800 & $0$  &$1 \;10^{-15}$           & $2 \;10^{-16}$      & $2 \;10^{-12}$   & $2 \;10^{-16}$   & $2 \;10^{-12}$  \\\hline
\end{tabular}\\
\bigskip
\caption{\it Relative errors for Clement matrices}
\label{tabelaRelClem}
\end{table}

\begin{table}[htbp]
\centering
\begin{tabular}{l||c|c|c|}
\multicolumn{1}{c}{$n$} & \multicolumn{1}{c}{$BGT$}         & \multicolumn{1}{c}{\textsl{complex dqds}} & \multicolumn{1}{c}{\textsl{tridqds}}  \\\hline \hline
100 & $4.2$     & $0.06$    & $0.06$       \\
200 & $12.6$    &$0.12$   & $0.13$      \\
400 & $42.2$      &$0.45$   & $0.50$        \\
800 & $174.2$      &$1.8$    & $1.8$        \\\hline
\end{tabular}\\
\bigskip
\caption{\it CPU times in seconds for Clement matrices}
\label{ClementAcc}
\end{table}

%%%%%%%%%%%%%%%%%%%%%%%%%%%%%%%%%%%%%%%%%%%%%%%%%%%%%%%%%%%%%%%%
\bigskip

\noindent\textbf{Graded matrix}

\bigskip

This matrix $C$ was created in $\Delta T$ form with
$T=\tridiag(\boldsymbol{b},\boldsymbol{a},\boldsymbol{c})$,
\begin{align*}
a_j&=b_j=c_j=3^{-(j-1)}, \quad j=1,\ldots,n-1; \quad a_n=3^{-(n-1)},
\end{align*}
and $\Delta=\diag(\boldsymbol{\delta})$, $\delta_{j}=(-1)^{\lfloor (j+1)/2 \rfloor}$, $j=1,\ldots,n$. The result is a balanced matrix with eigenvalues of different magnitude.

\begin{figure*}[htbp]
\centering{
   \includegraphics[width=9cm]{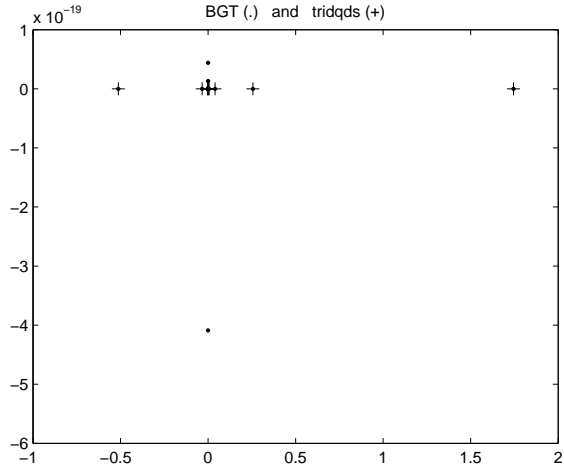}
   \caption{\it Eigenvalues of the graded matrix with $n=100$}
   \label{PictureGraded}
    }
\end{figure*}

Figure \ref{PictureGraded} shows the approximated eigenvalues for $n=100$. Table \ref{TableGraded} reports the CPU times.

\begin{table}[htbp]
\centering
\begin{tabular}{l||c|c|c|}
\multicolumn{1}{c}{$n$} & \multicolumn{1}{c}{$BGT$}         & \multicolumn{1}{c}{\textsl{complex dqds}} & \multicolumn{1}{c}{\textsl{tridqds}}  \\\hline \hline
50 & $0.31$    & $0.02$    & $0.02$       \\
100 &$0.67$    &$0.06$   & $0.04$      \\
200 &$2.12$    &$0.14$   & $0.06$        \\
400 & $...$    &$0.45$    & $0.35$        \\\hline
\end{tabular}\\
\bigskip
\caption{\it CPU times in seconds for the graded matrices}
\label{TableGraded}
\end{table}
$BGT$ code reported the message ``Exceed maximum number of operations" for $n=400$.
The performance of all methods for the flipped matrix is practically the same.

%%%%%%%%%%%%%%%%%%%%%%%%%%%%%%%%%%%%%%%%%%%%%%%%%%%%%%%%%%%%%%%%%%%%%%%%%%%%%%%%%%%
\pagebreak
\noindent\textbf{Matrix with clusters}

\bigskip

Matrix Test 5 in \cite{Tisseur},
\begin{align*}
C=D^{-1} \tridiag(\boldsymbol{1}, \boldsymbol{\alpha},\boldsymbol{1}), \;\;
D=\diag(\boldsymbol{\beta}), \quad \boldsymbol{\alpha},\boldsymbol{\beta} \in \matR^n   \;\; \; \\
\alpha_k=10^{5(-1)^k}\cdot(-1)^{\lfloor k/4\rfloor}, \;\; \beta_k=(-1)^{\lfloor k/3 \rfloor}, \quad k=1,\ldots,n,
\end{align*}
seems to be a challenging test matrix. It was designed to have large, tight clusters of eigenvalues around $10^{-5}$,
$-10^5$ and $10^5$. Half the spectrum is around $10^{-5}$ and the rest is divided unevenly between $-10^5$ and $10^5$.
The diagonal alternates between entries of absolute value $10^5$ and $10^{-5}$ and so, for \textsl{dqds} codes, there is a lot of rearranging
to do. When $n \geq 100$ it is not clear what is meant by accuracy.

All three codes obtain the correct number of eigenvalues in each cluster and the diameters of the clusters are all about $10^{-5}$.
The striking feature is the time taken. See Table \ref{Tisseur5Times}.

\begin{table}[htbp]
\centering
\begin{tabular}{l||c|c|c|}
\multicolumn{1}{c}{$n$} & \multicolumn{1}{c}{$BGT$}         & \multicolumn{1}{c}{\textsl{complex dqds}} & \multicolumn{1}{c}{\textsl{tridqds}}  \\\hline \hline
50  &$1.2$     & $0.03$    & $0.01$       \\
100 &$4.5$     & $0.05$    & $0.03$       \\
200 &$20.1$    &$0.14$     & $0.08$      \\
400 &$85.0$    &$0.61$     & $0.13$        \\\hline
\end{tabular}\\
\bigskip
\caption{\it CPU times in seconds for Test 5 matrix}
\label{Tisseur5Times}
\end{table}

\pagebreak
%%%%%%%%%%%%%%%%%%%%%%%%%%%%%%%%%%%%%%%%%%%%%%%%%%%%%%%%%%%%%%%%%%%%%%%%%%%%%%%%%%%
\bigskip
 \noindent\textbf{Matrix Test 4 }

\bigskip

\begin{figure*}[htbp]
\centering{
   \includegraphics[width=9cm]{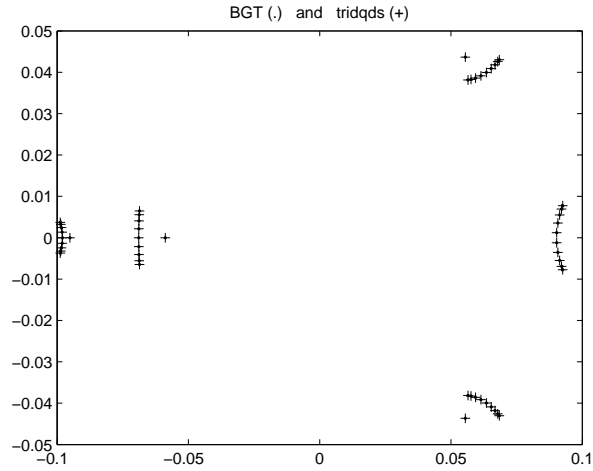}
   \caption{\it Eigenvalues of Test 4 matrix with $n=50$}
   \label{FigTeste4}
    }
\end{figure*}

For matrix Test 4 in \cite{Tisseur},
\begin{align*}
C=D^{-1} \tridiag(\boldsymbol{1}, \boldsymbol{\alpha},\boldsymbol{1}), \;\;
D=\diag(\boldsymbol{\beta}), \quad \boldsymbol{\alpha},\boldsymbol{\beta} \in \matR^n   \;\; \; \\
\alpha_k=(-1)^k, \;\; \beta_k=20\cdot(-1)^{\lfloor k/5 \rfloor}, \quad k=1,\ldots,n,
\end{align*}
the performance of the three codes is shown in Table \ref{TableTest4}. Figure \ref{FigTeste4} shows the eigenvalues of this matrix for $n=50$.

\begin{table}[htbp]
\centering
\begin{tabular}{l||c|c|c|}
\multicolumn{1}{c}{$n$} & \multicolumn{1}{c}{$BGT$}         & \multicolumn{1}{c}{\textsl{complex dqds}} & \multicolumn{1}{c}{\textsl{tridqds}}  \\\hline \hline
50  &$1.0$     & $0.05$    & $0.03$       \\
100 &$4.3$     & $0.08$    & $0.03$       \\
200 &$17.8$    &$0.27$     & $0.09$      \\
400 &$78.6$    &$1.0$     & $0.44$        \\
800 & $342.6$     &$5.5$      & $2.4$        \\\hline
\end{tabular}\\
\bigskip
\caption{\it CPU times in seconds for Test 4 matrix}
\label{TableTest4}
\end{table}

%%%%%%%%%%%%%%%%%%%%%%%%%%%%%%%%%%%%%%%%%%%%%%%%%%%%%%%%%%%%%%%%%%%%%%%%%%%%%%%%%%%
\bigskip
\noindent\textbf{Liu matrix}

\bigskip

Z. A. Liu \cite{Liu} devised an algorithm to obtain one-point spectrum unreduced tridiagonal matrices of
arbitrary dimension $n \times n$.
These matrices have only one eigenvalue, zero with multiplicity $n$, and the Jordan form consists of one Jordan block.
Our code \textsl{tridqds} computes this eigenvalue exactly (and also the generalized eigenvectors) using the
following method which is part of the prologue.

The best place to start looking for eigenvalues of a tridiagonal matrix
$C=\tridiag(\boldsymbol{a},\boldsymbol{b},\boldsymbol{c})$ is at the
arithmetic mean which we know $(\mu=\trace(C)/n)$. Before converting to $J-$form
and factoring, we check whether $\mu$ is an eigenvalue by using the 3-term recurrence
to solve
\[
(\mu I - C)\boldsymbol{x}=\boldsymbol{e}_np_n(\mu)/\prod_{i=1}^{n-1}c_i.
\]
Here
\begin{align*}
x_1&=1, \quad x_2=(\mu-a_2)/c_1,\\
x_{j+1}&=\frac{1}{c_j}\left[(\mu-a_j)x_j-b_{j-1}x_{j-1}\right], \quad j=2,\ldots,n-1,
\end{align*}
and
\[
\upsilon:=(\mu-a_n)x_n-b_{n-1}x_{n-1}\left(=p_n(\mu)/\prod_{i=1}^{n-1}c_i\right).
\]
If, by chance, $\upsilon$ vanishes, or is negligible compared to $\norm{\boldsymbol{x}}$,
then $\mu$ is an eigenvalue (to working accuracy) and $\boldsymbol{x}$ is an eigenvector.
To check its multiplicity we differentiate with respect to $\mu$ and solve
\[
(\mu I-C)\boldsymbol{y}=\boldsymbol{x}
\]
with $y_1=0$, $y_2=1=x_2^\prime$ $(=x_1)$. If
\[
\upsilon^\prime=p_n^\prime(\mu)/\prod_{i=1}^{n-1}c_i:=(\mu-a_n)y_n-b_{n-1}y_{n-1}+x_n
\]
vanishes, or is negligible w.r.t. $\norm{\boldsymbol{y}}$, then we continue the same way until the system
is inconsistent or there are $n$ generalized eigenvectors.

Usually $\upsilon \neq 0$ and the calculation appears to have been a waste.
This is not quite correct. In exact arithmetic, triangular factorization of
$\mu I-C$ or $\mu I-J$, where $J=\Delta C \Delta^{-1}$, will break down if, and only if,
$x_j$ vanishes for $1<j<n$. So our code examines $\min_j |x_j|$ and if it is too small
w.r.t. its neighbors and w.r.t. $\norm{\boldsymbol{x}}$ then we do not choose $\mu$ as our
initial shift. Otherwise we do obtain initial $L$ and $U$ from $J-\mu I=LU$.

%%%%%%%%%%%%%%%%%%%%%%%%%%%%%%%%%%%%%%%%%%%%%%%%%%%%%%%%%%%%%%%%%%%%%%%%%%%%%%%%%%%
\section{Conclusions and future work}\label{Section8}
We conclude that, working together, a single \textsl{dqds} transform with real shifts and
our \textsl{tridqds} transform with complex conjugate pairs of shifts constitute the right
tool for computing the eigenvalues of real tridiagonal matrices.

However there is far more work to be done for the following reasons. In a previous paper
we discovered that, surprisingly often, eigenvalues are determined to,
not high, but adequate relative accuracy; tiny relative changes $\eta$ in the parameters
that define the matrix produce relative changes in the eigenvalue of the order of $10^3\eta$
or $10^4\eta$. This is good news. We cannot tell in advance when this occurs. In our opinion a relative condition
number should be returned with each eigenvalue. This requires an approximation to the row and column eigenvectors,
whether or not the user needs them.

We envision software that computes an initial approximation to each eigenvalue and then invokes
a generalized Rayleigh quotient iteration to both compute eigenvectors and obtain a refined eigenvalue
approximation, along with the smallest residual norms that could be achieved.
Then the relative condition number can be formed. See \cite{Cferreira2}.

Another practical feature is to scan the initial matrix to extract Gersgorin disks and
a tight box in the complex plane that contains the spectrum. Matrices from industrial sources
frequently permit ``localization" of the eigenvectors belonging to certain parts of the spectrum.
One consequence is that the relevant eigenvectors, and eigenvalues, may be obtained from
small submatrices.

It is also important to scale and normalize the initial matrix and make use of the splitting
that occurs with big matrices. Currently our shift strategy is quite straightforward and it is both
difficult and worthwhile to improve it. There are plenty of challenges to be met before software
for this real tridiagonal problem can be installed in packages such as \textsc{LAPACK}, not to mention parallel
computation and \textsc{scaLAPACK}.

%%%%%%%%%%%%%%%%%%%%%%%%%%%%%%%%%%%%%%%%%%%%%%%%%%%%%%%%%%%%%%%%%%%%%%%%%%%%%%%%%%%

\appendix
\section{Implementation details of minor step \textit{i}}
\label{ApendiceA}
In this appendix we show how the calculations involved
in each minor step of \textsl{tridqds} can be organized.

For each minor step $i, \; i=2,\ldots,n-3$, consider $F$ and $G$ as
in (\ref{BeginStep}). Denote the $2 \times 1$ bulge in $F$,
indicated with plus signs, by $\begin{bmatrix} x_l &y_l
\end{bmatrix}^T$. And denote the entries $G_{i,i}$, $G_{i+1,i}$ and $G_{i+2,i}$, indicated with $*,+,+$, by
$\begin{bmatrix} x_r &y_r & z_r
\end{bmatrix}^T$. Subscripts $l$ and $r$ derive from ``left" and
``right", respectively. This way we have
\begin{align}
 F= \begin{bmatrix}
        \ddots&         &       &          &         &    \\
        \ddots&        1&      &          &         &    \\
               &\hat{l}_{i-1}   &u_i       &1        &         &    \\
               &x_l       &       &u_{i+1}   &1        &    \\
               &y_l       &       &          &u_{i+2}  &\ddots   \\
               &        &       &          &         &\ddots
  \end{bmatrix},
  \qquad
  G= \begin{bmatrix}
          \ddots & \ddots&       &        &        &    \\
                 &\hat{u}_{i-1}  &1       &        &       &    \\
                 &               &x_r       &        &       &    \\
                 &               &y_r       &1       &       &    \\
                 &               &z_r       &l_{i+1} &1      &    \\
                 &               &        &        &\ddots &\ddots
 \end{bmatrix}\nonumber
\end{align}
and the minor step $i$ can be accomplished using only these auxiliary variables.

\medskip

\noindent $\underline{\textbf{Minor step \textit{i}}}$

\smallskip

\begin{itemize}
 \item[\textbf{a)}]
     \begin{itemize}
        \item[$\bullet$] Matrices $Z_i^{-1}$ and $Z_i$

         \begin{align*}
          Z_i^{-1}= \begin{bmatrix}
                     \ddots &           &       &       &     &\\
                            &          1&       &      &       &        \\
                            &           &   \frac{1}{u_i}      & -\frac{1}{u_i}               &       &       \\
                            &           &               0&             1        &       &         \\
                            &                     &               &       &1      &      \\
                            &                     &               &       &       &\ddots \\
                      \end{bmatrix},
                      &\qquad
         Z_i^{-1}= \begin{bmatrix}
                     \ddots &           &       &       &     &\\
                            &          1&       &      &       &        \\
                            &           &       u_i      & 1               &       &       \\
                            &           &               0&             1        &       &         \\
                            &                     &               &       &1      &      \\
                            &                     &               &       &       &\ddots \\
                      \end{bmatrix}
         \end{align*}

         \item[$\bullet$] The effect of $Z_i^{-1}$ and the effect of $Z_i$

         \begin{align*}
            FZ_i^{-1}= \begin{bmatrix}
                \ddots&         &       &          &         &    \\
                \ddots&        1&      &          &         &    \\
                  &\hat{l}_{i-1}   &1         &0        &         &    \\
                  &x_l      &       &u_{i+1}   &1        &    \\
                  &y_l       &       &          &u_{i+2}  &\ddots   \\
                  &        &       &          &         &\ddots
           \end{bmatrix},
           &\;
           Z_iG=\begin{bmatrix}
             \ddots & \ddots&       &        &        &    \\
                   &\hat{u}_{i-1}  &1       &        &       &    \\
                   &               &x_r       &1        &       &    \\
                   &               &y_r       &1       &       &    \\
                   &               &z_r       &l_{i+1} &1      &    \\
                   &               &        &        &\ddots &\ddots
                \end{bmatrix}
         \end{align*}

where
$$
  x_r \longleftarrow x_r*u_i+y_r\\
$$
\end{itemize}
%%%%%%%%%%%%%%%%%%%%%%%%%%%%%%%%%%%%%%%%%%%%%%%%%%%%%%%%
%%%%%%%%%%%%%%%%%%%%%%%%%%%%%%%%%%%%%%%%%%%%%%%%%%%%%%%%

 \item[\textbf{b)}]
  \begin{itemize}

   \item[$\bullet$]Matrices $\mathcal{L}_i^{-1}$ and $\mathcal{L}_i$

    \begin{align*}
    \mathcal{L}_i^{-1}= \begin{bmatrix}
        \ddots&         &       &          &         &    \\
              & \ddots      &      &          &         &    \\
               &       &1       &         &         &    \\
               &       &x_l       &1   &        &    \\
               &       &y_l       &          &1  &   \\
               &        &       &          &         &\ddots
    \end{bmatrix},
    \qquad
    \mathcal{L}_i= \begin{bmatrix}
        \ddots&         &       &          &         &    \\
              & \ddots      &      &          &         &    \\
               &       &1       &         &         &    \\
               &       &-x_l       &1   &        &    \\
               &       &-y_l       &          &1  &   \\
               &        &       &          &         &\ddots
    \end{bmatrix}
   \end{align*}

  where
  $$\begin{array}{l}
   x_l  \longleftarrow  -x_l/\hat{l}_{i-1} \\
   y_l  \longleftarrow  -y_l/\hat{l}_{i-1} \\
 \end{array}$$

   \item[$\bullet$]The effect of $\mathcal{L}_i^{-1}$

   $$
   \mathcal{L}_i^{-1}F= \begin{bmatrix}
        \ddots&         &       &          &         &    \\
        \ddots&        1&      &          &         &    \\
               &\hat{l}_{i-1}   &1         &        &         &    \\
               &        &x_l       &u_{i+1}   &1        &    \\
               &        &y_l      &          &u_{i+2}  &\ddots   \\
               &        &       &          &         &\ddots
  \end{bmatrix}
  $$

   \item[$\bullet$]The effect of $\mathcal{L}_i$

   $$
   G\mathcal{L}_i= \begin{bmatrix}
          \ddots & \ddots&       &        &        &    \\
                 &\hat{u}_i       &1        &       &    \\
                 &x_r       &1       &       &    \\
                 &y_r       &l_{i+1} &1      &    \\
                 &z_r        &        &l_{i+2} &1\\
                 &               &        &        &\ddots &\ddots
 \end{bmatrix}
  $$

  where
  $$\begin{array}{l}
  \hat{u}_i \longleftarrow x_r-x_l\\
  x_r \longleftarrow y_r-x_l\\
  y_r \longleftarrow z_r-y_l-x_l*l_{i+1}\\
  z_r \longleftarrow -y_l*l_{i+2}\\
  \end{array}$$

\end{itemize}
%%%%%%%%%%%%%%%%%%%%%%%%%%%%%%%%%%%%%%%%%%%%
%%%%%%%%%%%%%%%%%%%%%%%%%%%%%%%%%%%%%%%%%%%%%

 \item[\textbf{c)}]
  \begin{itemize}

  \item[$\bullet$]Matrices $Y_i^{-1}$ and $Y_i$
   \begin{align*}
   Y_i^{-1}= \begin{bmatrix}
          \ddots &           &   &        &        &    \\
                 &1          &        &       &    \\
                 &x_r        &1       &       &    \\
                 &y_r        &  &1      &    \\
                 &z_r        &        &  &1\\
                 &               &        &        & &\ddots
  \end{bmatrix},
  \qquad
  Y_i= \begin{bmatrix}
          \ddots &           &   &        &        &    \\
                 &1          &        &       &    \\
                 &-x_r        &1       &       &    \\
                 &-y_r        &  &1      &    \\
                 &-z_r        &        &  &1\\
                 &               &        &        & &\ddots
  \end{bmatrix}
  \end{align*}
 where
 $$\begin{array}{l}
 x_r \longleftarrow x_r/\hat{u}_i\\
 y_r \longleftarrow y_r/\hat{u}_i\\
 z_r \longleftarrow z_r/\hat{u}_i\\
 \end{array}$$

\item[$\bullet$] The effect of $Y_i^{-1}$
 $$
 FY_i^{-1}= \begin{bmatrix}
        \ddots&         &       &          &         &    \\
        \ddots&        1&      &          &         &    \\
               &\hat{l}_i   &u_{i+1}       &1        &         &    \\
               &x_l      &       &u_{i+2}   &1        &    \\
               &y_l       &       &          &u_{i+3}  &\ddots   \\
               &        &       &          &         &\ddots
  \end{bmatrix}
 $$
 where
  $$\begin{array}{l}
  \hat{l}_i \longleftarrow x_l+y_r+x_r*u_{i+1}\\
   x_l  \longleftarrow  y_l+z_r+ y_r*u_{i+2}\\
   y_l  \longleftarrow  z_r*u_{i+3} \\
 \end{array}$$

 \item[$\bullet$]The effect of $Y_i$
 $$
 Y_iG= \begin{bmatrix}
          \ddots & \ddots&       &        &        &    \\
                 &\hat{u}_i      &1       &        &       &    \\
                 &               &x_r       &        &       &    \\
                 &               &y_r       &1       &       &    \\
                 &               &z_r       &l_{i+2} &1      &    \\
                 &               &        &        &\ddots &\ddots
 \end{bmatrix}
 $$
 where
 $$\begin{array}{l}
 x_r \longleftarrow 1-x_r\\
 y_r \longleftarrow l_{i+1}-y_r\\
 z_r \longleftarrow -z_r\\
 \end{array}$$
\end{itemize}
\end{itemize}

%%%%%%%%%%%%%%%%%%%%%%%%%%%%%%%%%%%%%%%%%%%%%%%%%%%%%%%%%%%%%%%%%%%
\section{\textit{tridqds} algorithm}
\label{ApendiceB}

\bigskip

$$
\begin{array}{rll}
\textsl{tridqds}(\sigma, \bar{\sigma}):& \\\\
             &\textsf{\% \text{ step } 1} \\  % STEP 1
             & x_r = 1; \;\; y_r = l_1; \;\; z_r = 0 \\

             & \textsf{\%  \text{ the effect of $Z_1$}} \\
             & x_r = x_r*u_1+y_r \\

             & \textsf{\% \text{ the matrix $\mathcal{L}_1^{-1}$}}\\
             & x_l=(u_1+l_1)^2 + u_2l_1 - 2 (\Re\sigma)(u_1+l_1)+|\sigma|^2\\
             & y_l = -u_2l_1u_3l_2/x_l \\
             & x_l = -u_2l_1(u_1+l_1 + u_2+l_2 - 2 (\Re\sigma))/x_l\\

             & \textsf{\% \text{ the effect of $\mathcal{L}_1$}}\\
             & \hat{u}_1=x_r-x_l;\\
             & x_r=y_r-x_l; \;\; y_r=z_r-y_l-x_l*l_2; \;\; z_r=-y_l*l_3\\

             & \textsf{\% \text{ the matrix $Y_1^{-1}$}}\\
             & x_r = x_r/\hat{u}_1; \;\; y_r = y_r/\hat{u}_1; \;\; z_r = z_r/\hat{u}_1\\

             & \textsf{\% \text{ the effect of $Y_1^{-1}$}}\\
             & \hat{l}_1 = x_l+y_r+x_r*u_2\\
             & x_l = y_l+z_r+ y_r*u_3; \;\; y_l = z_r*u_4\\

             & \textsf{\% \text{ the effect of $Y_1$}}\\
             & x_r = 1-x_r; \;\; y_r = l_2-y_r; \;\; z_r =-z_r\\\\

\end{array}
$$

$$
\begin{array}{rll}
\hspace{3.5cm}&\\
             & \textbf{for } i=2, \ldots, n-3 \\   %FOR
             & \hspace{0.6 cm} \textsf{\% \text{ the effect of $Z_i$}} \\
             & \hspace{0.6 cm} x_r =x_r*u_i+y_r   \\

             & \hspace{0.6 cm} \textsf{\% \text{ the matrix $\mathcal{L}_i^{-1}$}}\\
             & \hspace{0.6 cm} x_l= -x_l/\hat{l}_{i-1} ; \;\; y_l=-y_l/\hat{l}_{i-1};\;\;\\

             & \hspace{0.6 cm} \textsf{\% \text{ the effect of $\mathcal{L}_i$}}\\

             & \hspace{0.6 cm} \hat{u}_i=x_r-x_l;\\
             & \hspace{0.6 cm} x_r=y_r-x_l; \;\; y_r=z_r-y_l-x_l*l_{i+1}; \;\; z_r=-y_l*l_{i+2}\\

             & \hspace{0.6 cm} \textsf{\% \text{ the matrix $Y_i^{-1}$}}\\
             & \hspace{0.6 cm} x_r = x_r/\hat{u}_i; \;\; y_r = y_r/\hat{u}_i; \;\; z_r = z_r/\hat{u}_i\\

             & \hspace{0.6 cm} \textsf{\% \text{ the effect of $Y_i^{-1}$}}\\
             & \hspace{0.6 cm} \hat{l}_i = x_l+y_r+x_r*u_{i+1}\\
             & \hspace{0.6 cm} x_l = y_l+z_r+ y_r*u_{i+2}; \;\; y_l = z_r*u_{i+3}\\

             & \hspace{0.6 cm} \textsf{\% \text{ the effect of $Y_i$}}\\
             & \hspace{0.6 cm} x_r = 1-x_r; \;\; y_r = l_{i+1}-y_r; \;\; z_r =-z_r\\

             & \textbf{end for}\\\\     %END FOR
\end{array}
$$
$$
\begin{array}{rll}
\hspace{3.5cm}&\\
             &\textsf{\% \text{ step } n-2} \\  % STEP N-2
             & \textsf{\% \text{ the effect of $Z_{n-2}$}} \\
             & x_r = x_r*u_{n-2}+y_r   \\

             & \textsf{\% \text{ the matrix $\mathcal{L}_{n-2}^{-1}$}}\\
             & x_l= -x_l/\hat{l}_{n-3}  ; \;\; y_l = -y_l/\hat{l}_{n-3};\;\;\\

             & \textsf{\% \text{ the effect of $\mathcal{L}_{n-2}$}}\\
             & \hat{u}_{n-2}=x_r-x_l;\\
             & x_r=y_r-x_l; \;\; y_r= z_r-y_l-x_l*l_{n-1} \\ %%%%%

             & \textsf{\% \text{ the matrix $Y_{n-2}^{-1}$}} \hspace{4.5cm} \\
             & x_r = x_r/\hat{u}_{n-2} ; \;\; y_r = y_r/\hat{u}_{n-2} \\

             & \textsf{\% \text{ the effect of $Y_{n-2}^{-1}$}}\\
             & \hat{l}_{n-2} = x_l+y_r+x_r*u_{n-1} \\
             & x_l = y_l+ y_r*u_{n}  \\

             & \textsf{\% \text{ the effect of $Y_{n-2}$}}\\
             & x_r = 1-x_r; \;\; y_r = l_{n-1}-y_r  \\\\
\end{array}
$$

$$
\begin{array}{rll}
             &\textsf{\%\text{ step } n-1 } \\             % STEP N-1
             & \textsf{\% \text{ the effect of $Z_{n-1}$}} \\
             & x_r = x_r*u_{n-1}+y_r    \\

             & \textsf{\% \text{ the matrix $\mathcal{L}_{n-1}^{-1}$}}\\
             & x_l= -x_l/\hat{l}_{n-2} \\

             & \textsf{\% \text{ the effect of $\mathcal{L}_{n-1}$}}\\
             & \hat{u}_{n-1}=x_r-x_l;\\
             & x_r=y_r-x_l  \\ %%%%%

             & \textsf{\% \text{ the matrix $Y_{n-1}^{-1}$}}\\
             & x_r = x_r/\hat{u}_{n-1} \\

             & \textsf{\% \text{ the effect of $Y_{n-1}^{-1}$}}\\
             & \hat{l}_{n-1} = x_l+x_r*u_n  \\

             & \textsf{\% \text{ the effect of $Y_{n-1}$}}\\
             & x_r = 1-x_r \\\\
\end{array}
$$
$$
\begin{array}{rll}
             &\textsf{\textsf{\%\text{ step } n}} \\     % STEP N
             & \% \textsf{ the effect of $Z_n$} \\
             & x_r = x_r*u_n    \\

             & \textsf{\% \text{ the matrix $\mathcal{L}_n^{-1}=I$}}\\

             & \textsf{\% \text{ the effect of $\mathcal{L}_n$}}\\
             & \hat{u}_n= x_r;\\

             & \textsf{\% \text{ the matrix $Y_n^{-1}=I$}}\\
\end{array}
$$
\medskip
\pagebreak

%%%%%%%%%%%%%%%%%%%%%%%%%%%%%%%%%%%%%%%%%%%%%%%%%%%%%%%%%%%%%%%%%%%%%%%%%%%%%%%%%%%

\end{document}

%% file: FiguraRelLRqdsTalkTese.tex
%%%%%%%%%%%%%%%%%%%%%%%%%%%%%%%%%%%%%%%%%%%%%%%%%%%%%%%%%%
%% Figura relativa \`{a} rela\c{c}\~{a}o dos m\'{e}todos LR e qds
%%%%%%%%%%%%%%%%%%%%%%%%%%%%%%%%%%%%%%%%%%%%%%%%%%%%%%%%%%

\setlength{\unitlength}{0.8 cm}

\begin{center}

\begin{picture}(13.3,4)

\put(-1,3.3) {\makebox(1,1)[b]{$J_1$}}

\put(0,3.5){\vector(1,0){3}}

\put(3,3.3) {\makebox(1,1)[b]{$J_2$}}

\put(4,3.5){\vector(1,0){3}}

\put(7,3.3) {\makebox(1,1)[b]{$J_3$}}

\put(8,3.5){\vector(1,0){3}}

\put(11,3.3) {\makebox(1,1)[b]{$J_4$}}

%%%%%%%%%%%%%%%%%%%%%%%%%%%%%%%%%%%%%%%%%%%%%%%%%%%

\put(0.5,3.6) {\makebox(2,1)[b]{$\text{LR}(0)$}}

\put(4.5,3.6) {\makebox(2,1)[b]{$\text{LR}(\sigma)$}}

\put(8.5,3.6) {\makebox(2,1)[b]{$\text{LR}(\bar{\sigma})$}}

%%%%%%%%%%%%%%%%%%%%%%%%%%%%%%%%%%%%%%%%%%%%%%%%%%%

\put(-0.2,3){\vector(1,-1){1.5}}

\put(1.8,1.5){\vector(1,1){1.5}}

\put(3.8,3){\vector(1,-1){1.5}}

\put(5.8,1.5){\vector(1,1){1.5}}

\put(7.8,3){\vector(1,-1){1.5}}

\put(9.8,1.5){\vector(1,1){1.5}}

\put(11.8,3){\vector(1,-1){1.5}}
%%%%%%%%%%%%%%%%%%%%%%%%%%%%%%%%%%%%%%%%%%%%%%%%%%%%%

\put(0.5,0.7){\makebox(2,1)[b]{$L_1,U_1$}}

\put(2.3,1){\vector(1,0){2.4}}

\put(4.5,0.7){\makebox(2,1)[b]{$L_2,U_2$}}

\put(6.3,1){\vector(1,0){2.4}}

\put(8.5,0.7){\makebox(2,1)[b]{$L_3,U_3$}}

\put(10.3,1){\vector(1,0){2.4}}

\put(12.5,0.7){\makebox(2,1)[b]{$L_4,U_4$}}

%%%%%%%%%%%%%%%%%%%%%%%%%%%%%%%%%%%%%%%%%%%%%%

\put(2.5,0.2){\makebox(2,1)[b]{$\text{dqds}(\sigma)$}}

\put(6.5,0.2){\makebox(2,1)[b]{$\text{dqds}(\overline{\sigma}-\sigma)$}}

\put(10.5,0.2){\makebox(2,1)[b]{$\text{dqds}(-\overline{\sigma})$}}

\put(4.2,2.3){\makebox(1,1)[b]{$-\sigma$}}

\put(5.8,2.3){\makebox(1,1)[b]{$\sigma$}}

\put(8.2,2.3){\makebox(1,1)[b]{$-\overline{\sigma}$}}

\put(9.8,2.3){\makebox(1,1)[b]{$\overline{\sigma}$}}

\end{picture}

\end{center}

%%%%%%%%%%%%%%%%%%%%%%%%%%%%%%%%%%%%%%%%%%%%%%%%%%%%%%%%%%%%%

%% file: CommutativeDiagrame.tex
%COMMUTATIVE DIAGRAM FOR dqds ALGORITHM

\newcommand{\textA}{\operatorname{\scriptsize{\begin{array}{l} \text{change each} \\
                                                         l_k \text{ by 1 ulp}\\
                                                         u_k \text{ by 3 ulps}\\
                                                     \end{array}}}}
\newcommand{\textB}{\operatorname{\scriptsize{\begin{array}{l} \text{change each} \\
                                                         \breve{l}_k,\; \breve{u}_k  \text{ by 2 ulps}\\
                                                     \end{array}}}}
\newcommand{\dqds}{\operatorname{\hspace{0.5cm} dqds \hspace{0.5cm}}}
\newcommand{\computed}{\operatorname{computed}}
\newcommand{\exact}{\operatorname{exact}}

\begin{equation*}
\begin{CD}
 {L,U}             @>{\dqds}>{\computed}>    {\widehat{L},\widehat{U}}\\
 @V{\textA}VV                                   @AA{\textB}A \\
 {\widetilde{L},\widetilde{U}} @>{\dqds}>{\exact}>   {\breve{L},\breve{U}}
\end{CD}
\end{equation*}